\documentclass{amsart}

\begin{document}

\title[Hasse polynomial]{Lie Invariant  Frobenius lifts and \\ deformation of the Hasse polynomial}

\def \h{\hat{\ }}
\def \cO{\mathcal O}
\def \ra{\rightarrow}
\def \bZ{{\mathbb Z}}
\def \cP{\mathcal V}
\def \cH{{\mathcal H}}
\def \cB{{\mathcal B}}
\def \d{\delta}
\def \cC{{\mathcal C}}
\def \jor{\text{jor}}

\newtheorem{THM}{{\!}}[section]
\newtheorem{THMX}{{\!}}
\renewcommand{\theTHMX}{}
\newtheorem{theorem}{Theorem}[section]
\newtheorem{corollary}[theorem]{Corollary}
\newtheorem{lemma}[theorem]{Lemma}
\newtheorem{proposition}[theorem]{Proposition}
\newtheorem{thm}[theorem]{Theorem}
\theoremstyle{definition}
\newtheorem{definition}[theorem]{Definition}
\theoremstyle{remark}
\newtheorem{remark}[theorem]{Remark}
\newtheorem{example}[theorem]{Examples}
\numberwithin{equation}{section}
\newtheorem{conjecture}[theorem]{\bf Conjecture}
\address{Department of Mathematics and Statistics\\University of New Mexico \\
  Albuquerque, NM 87131, USA\\ 
   Department of Mathematics} 
\subjclass[2010]{Primary: 11G07, Secondary: 11F33}
\maketitle

\bigskip

\medskip
\centerline{\bf Alexandru Buium}
\bigskip

\begin{abstract} 
We show  that the $p$-adic completion of any affine 
 elliptic curve with ordinary reduction
 possesses Frobenius lifts
whose ``normalized" action on $1$-forms 
preserves mod $p$
 the space of invariant $1$-forms.
We next  show that, after removing the $2$-torsion sections, 
the above situation can be ``infinitesimally deformed" in the sense that
the above mod $p$ result has a mod $p^2$ analogue.
 While the ``eigenvalues" mod $p$ are given by 
the reciprocal of the Hasse polynomial, 
the ``eigenvalues" mod $p^2$ are given by an appropriate  $\d$-modular function whose reciprocal is a $p$-adic deformation 
 of  the Hasse polynomial. 
\end{abstract}

\bigskip
\bigskip
\bigskip

\section{Introduction}
Consider a complex affine elliptic curve, its smooth projective compactification,  a basis $\omega$ for the global (equivalently, translation invariant) $1$-forms on the compactification, and a basis $\epsilon$ for the  translation invariant vector fields.
For an algebraic vector field $\d$ on our affine elliptic curve the following are well known to be equivalent: 

1) $\d$ extends to  the compactification;

2)  the Lie derivative along $\d$
 kills $\omega$;

3) $\d$ commutes with $\epsilon$; 
 
 4) $\d$ is a constant times $\epsilon$.

We would like to explore an arithmetic analogue of this picture. As such, the present paper is part of a more general project devoted to investigating arithmetic analogues of  classical differential equations; cf., especially, \cite{char, difmod, book, BYM, euler, canonical}. Indeed, conditions 1)-4) above appear, for instance,  as the 
``linearization condition" for certain completely integrable systems such as 
 the Euler equations for the rigid body; and an  arithmetic analogue of the Euler equations was introduced and studied in \cite{euler, canonical}. The present paper is mainly motivated by  questions raised in  \cite{euler, canonical} but it is written so as to be independent of all of the above cited papers. 
 
 For our arithmetic setting let us take, as our ground ring, a complete discrete valuation ring $R$ with maximal ideal generated by a prime $p$ and algebraically closed residue field $\overline{R}:=R/pR$; such an $R$ is isomorphic to the Witt ring on $\overline{R}$.  To simplify working with elliptic curves we will assume $p\neq 2,3$. As an analogue of vector fields on complex algebraic varieties one can take Frobenius lifts
 (i.e. endomorphisms lifting the characteristic $p$ Frobenius) on $p$-adic formal
schemes  over  $R$; cf. \cite{char, book} where this analogy is part of a more general strategy. Note at this point, however, that Frobenius lifts on $p$-adic completions of affine 
elliptic curves over $R$ rarely extend to the compactification of our elliptic curves; the existence of such an extension for an elliptic curve forces that curve to have complex multiplication. So the arithmetic analogue of condition 1) above is too restrictive. Also note that condition 4) above is meaningless in the arithmetic context.  It turns out, however,  that  (equivalent) arithmetic analogues of conditions 2) and 3) can be introduced. For simplicity let us explain the arithmetic analogue of condition 2) only.  Let $E$ be the {\it affine plane elliptic curve} over $R$ defined by
$$E=Spec\ R[x,y]/(y^2-(x^3+ax+b))$$
 and consider the $1$-form $\omega=\frac{dx}{y}$ on $E$; this form extends to an (automatically translation invariant) $1$-form $\omega$ on the smooth projective model  of $E$.
Let us say that a Frobenius lift $\phi$
on the $p$-adic completion 
$\widehat{E}$ of $E$
 is  {\it Lie invariant mod $p$} (with eigenvalue $\lambda\in R$) if the following holds:
\begin{equation}
\label{inver}
\frac{\phi^*}{p}\omega\equiv \lambda \omega\ \ \ \text{mod}\ \ \ p.\end{equation}
This concept played a role in \cite{BYM, euler, canonical}.
In \cite{foundations}, pp. 152, another concept of {\it invariance} for a Frobenius lift was considered; that concept will play no role in the present paper.
Our starting point is the following easy result  (which is a consequence of Theorem \ref{thm1}):

\begin{thm}
\label{easyy}
For any ordinary  $E$  there exists a Frobenius lift $\phi$
on $\widehat{E}$ such that $\phi$ is  Lie invariant mod $p$, with eigenvalue  $\lambda=H(a,b)^{-1}$, where $H$ is  the Hasse polynomial.
\end{thm}

A  variant of Theorem \ref{thm1}  is implicit in \cite{euler}. 
 One can view the congruence \ref{inver} as an analogue mod $p$ of condition 2) above; indeed the ``normalized" action of $\phi$ on forms (by which we will mean the action of $\frac{\phi^*}{p}$) can be viewed as  an analogue of the action on forms of the Lie derivative along a vector field; and  condition  \ref{inver} can be viewed as  an analogue of the condition that
 a Lie derivative annihilates  a form; cf. \cite{BYM, euler, canonical} where this analogy plays a role. 
 
 Now, following \cite{euler, canonical}, one can raise the question
whether, at least after replacing $E$ by an affine open set $E^*$ of it, we have that   \ref{inver} can be ``deformed infinitesimally" i.e., whether \ref{inver} can be lifted to a congruence,
\begin{equation}
\label{68}\frac{\phi^*}{p}\omega\equiv \lambda \omega\ \ \text{mod}\ \ p^2,
\end{equation}
where $\lambda$ is some element in $R$; 
we can view \ref{68} as an analogue mod $p^2$ of the condition 2) in the classical case. 
If \ref{68} holds we will  say that $\phi$ is {\it Lie invariant mod $p^2$} (with {\it eigenvalue} $\lambda$). 
Let $E^*$ be the complement in $E$ of the three $2$-torsion points and let $\widehat{E^*}$ be the $p$-adic completion of $E^*$. 
Our next main result (Theorem \ref{godem5}) will imply the following:

 \begin{thm}\label{comics}
 Assume $p\geq 13$ and $p\not\equiv 11$ mod $12$. 
 There is a finite subset  ${\mathcal F}\subset \overline{R}$ with the following property. 
 For any  $E$ with $j$-invariant mod $p$ not in ${\mathcal F}$
there is a Frobenius lift $\phi$ on $\widehat{E^*}$  such that $\phi$ is Lie invariant mod $p^2$. 
\end{thm}

It is not clear whether one should expect  Theorem \ref{comics} to hold with  $\widehat{E^*}$ replaced by $\widehat{E}$.
On the other hand one can ask if Theorem \ref{comics} holds for 
 other  primes; 
 also one can ask if ${\mathcal F}$ 
  can  be taken to be the non-ordinary (i.e., supersingular) locus;
  cf. Corollary \ref{godem4} for some special cases when these questions have positive answers.
 
Theorems \ref{thm1} and \ref{godem5} will be  relative versions of Theorems \ref{easyy} and \ref{comics}. In particular Theorem \ref{godem5} will show that, in Theorem \ref{comics}, the eigenvalue of $\phi$ will be given by $\Lambda(a,b)$ where $\Lambda$ will be a {\it singular $\d$-modular function of weak weight $1-p$} in the sense of \cite{difmod}. We will review these concepts in the body of the paper; we will also introduce in our paper some new concepts related to $\d$-modular forms (such as {\it quasi linearity}) and will prove that $\Lambda$ is  quasi linear. For the reader familiar with \cite{difmod} note that
the reciprocal of our $\d$-modular function $\Lambda$ will be congruent mod $p$ to the  Hasse polynomial $H$
and hence  can be viewed as a natural  ``$p$-adic deformation" of $H$. Now,
as well known, $H$ viewed as a modular function of weight $p-1$, has Fourier expansion congruent to $1$ mod $p$; cf \cite{Katzpadic}. 
On the other hand, in \cite{Barcau, book}, a unique $p$-adic deformation of $H$ to a $\d$-modular function of {\it weight} $\phi-1$, denoted by $f^{\partial}$, was proved to exist, 
with the property that the {\it $\d$-Fourier expansion} of $f^{\partial}$ is equal (rather than congruent) to $1$;
it would be interesting to understand whether there is a  connection  between the two deformations $\Lambda^{-1}$ and $f^{\partial}$ of $H$. 

  The strategy of our proof of Theorem \ref{comics} is as follows. 
  Given a Frobenius lift $\phi$ on $\widehat{E^*}$ it is easy to see that
  $\phi$ is Lie invariant mod $p^2$ if
  a certain rational function $Z(x)$ attached to $\phi$ is a solution mod $p^2$
  to a certain order $1$ linear differential equation with   coefficients  rational functions of $x$. We will look for solutions $Z(x)$ of a certain  special shape, cf. equation 
  \ref{gogol}; this allows one to reduce the problem of finding $Z(x)$  to showing  that
  a certain ``universal''   infinite  triangular system of algebraic (i.e., not differential)
   linear congruences mod $p$
  in infinitely many variables has a solution  ``with finite support''; cf. Proposition \ref{operastabbing}. 
  The above infinite system turns out to have a somewhat unexpected property:   solutions with finite support can be constructed by  solving a  ``truncated"  version of the  system in which only the first $\frac{p+7}{2}$ equations play a role; cf. Proposition \ref{stabcr}. 
  We conclude by showing 
  that this truncated system is solvable, at least if $p\not\equiv 11$ mod $12$, $p\geq 13$; cf. Proposition \ref{golem}. 
  The  strategy outlined above is sufficient to prove Theorem \ref{comics} above. A further layer of easy arguments, involving $\d$-modular functions, 
  is  necessary to prove the version that holds in families; cf. Theorem \ref{godem5}. 
  
  A similar (but much easier) argument can be employed to prove Theorem \ref{easyy}, at least for $E$ replaced by $E^*$; passing from $E^*$ to $E$ is possible due to a somewhat unexpected cancellation of a denominator in a key equation; cf. \ref{maduc}.

   If instead of the congruence mod $p^2$ we were looking at we want to explore congruences mod $p^m$ for arbitrary $m\geq 1$ our method leads to the consideration of a differential equation mod $p^m$ which is not linear anymore but, rather, resembles the differential equation satisfied by Weierstrass elliptic functions. In our differential equation  a quadratic polynomial in the derivative $\frac{d Z}{dx}$ of the unknown function $Z=Z(x)$ will be congruent to a cubic polynomial in $Z$; all these polynomials will have coefficients rational functions of $x$. On the other hand, in the classical Weierstrass equation for complex elliptic curves with uniformizing complex variable $w\in {\mathbb C}$, the square of the derivative $\frac{dx}{dw}$ of the unknown function $x=x(w)$  is a cubic polynomial in $x$ with constant coefficients. Note by the way the contrast between the independent variables in the two cases: in the classical case $\frac{d}{dw}$ is the derivative with respect to the uniformizing variable $w$ for the elliptic curve; in our case $\frac{d}{dx}$ is the derivative with respect to the coordinate $x$ on the elliptic curve. 
   Now the cubic polynomial $G$ becomes  quadratic in $Z$ when reduced mod $p^3$ and  becomes (together with its square root) linear in $Z$ when reduced mod $p^2$. In particular for $m=2$ our differential equation satisfied by $Z$ becomes linear. 
   The non-linear differential equations satisfied by $Z$ mod $p^m$ for $m\geq 3$ seem elusive at this time and this is why our paper concentrates on the cases $m=1,2$

  As explained above the problem of finding Frobenius lifts that are Lie invariant mod powers of $p$ boils down to finding  rational functions that are solutions mod powers of $p$ to certain special differential equations with coefficients rational functions. Such solutions rarely exist for general differential equations so the fact that, in our special case, such solutions exist is, in some sense, unexpected. (The existence of solutions in Theorem \ref{easyy} is more easily explained but for Theorem \ref{comics} the explanation is rather involved.)
  By the way, the theory of elliptic curves provides another striking example of a linear differential equation, with coefficients  rational functions,  satisfied mod $p$ by a rational function: indeed recall the result (going back to Deuring and Manin, cf. \cite{silverman}, pp. 140-143) that
  the Hasse polynomial $h(t)$ of the Legendre elliptic curve 
  $y^2=x(x-1)(x-t)$ satisfies the Picard-Fuchs 
  differential equation mod $p$ involving the differential operator $\frac{d}{dt}$. 
  One should note again the difference between such a statement and the paradigm of our proofs described above: in our paradigm
  the independent variable in the function $Z(x)$ and in the differential operator $\frac{d}{dx}$ is $x$; in the Deuring-Manin paradigm the independent variable in $h(t)$ and $\frac{d}{dt}$ is $t$.
  
  As final remark on Theorems \ref{thm1} and \ref{comics} note that 
  the reduction mod $p$ of our $E$ is a {\it hyperbolic curve} in the sense of Mochizuchi \cite{mochizuki}. But note that our $\widehat{E}$  is an aribitrary lift of its reduction mod $p$ and hence is generally  not equal to  {\it Mochizuki's canonical lift} of its reduction mod $p$. In particular Mochizuki's  theory \cite{mochizuki} (cf. also \cite{finotti}) of canonical Frobenius lifts on canonical lifts of ordinary hyperbolic curves does not apply to our situation.  Note however that a certain polynomial $W(x)$ that plays a role here (cf. \ref{grab}) also plays a role in \cite{finotti}, p. 25.

  All  the discussion above was devoted to elliptic curves.
 One can ask, at this point, whether  our results can be extended to other group schemes. 
 Let $G$ be a smooth group scheme over $R$. A Frobenius lift $\phi$ on the $p$-adic completion of an open subset $G^0\subset G$ is said to be  {\it Lie invariant mod $p$}  if \ref{inver} holds, where   $\omega$ 
 is  a column vector with entries a basis for the left invariant $1$-forms on $G$ and $\lambda$ is a square matrix with entries in $R$. (A similar definition can be made mod $p^m$.) 
 One is tempted to conjecture that any ordinary Abelian scheme $G$ has an open set whose $p$-adic completion possesses a 
    Frobenius lift that is Lie invariant mod $p$.  Cf. Conjecture \ref{milk}. (Theorem \ref{thm1} proves this  in the case of dimension $1$, when the  Frobenius lift is actually defined on the complement of the zero section.)  Also it is trivial to check that the multiplicative group $G={\mathbb G}_m=GL_1$,
 and hence any split torus, has a  Frobenius lift that is Lie invariant mod $p$ (defined on the whole of the group). In a separate paper \cite{alie2} we will show that {\it linear tori are  the only linear algebraic groups over a number field
 that admit Lie invariant Frobenius lifts for infinitely many primes}.
 
 The paper is organized as follows: in section 2 we  introduce the relevant concepts and notation; in section 3 we state our main results;
section 4 is devoted to the proofs of our results.

\medskip

{\bf Acknowledgments}. The present work was partially supported by the  IHES and  the Simons Foundation (award 311773). The author is grateful to Emma Previato for inspiring conversations.

\section{Main concepts}

\subsection{Frobenius lifts} Let  $p$ be a rational prime. 
Throughout the paper we assume $p\neq 2,3$.
 If $S$ is a ring a {\it Frobenius lift} 
on $S$  is a ring endomorphism $\phi:S\ra S$ whose reduction mod $p$ is the $p$-power 
Frobenius on $S/pS$. 
 Given a Frobenius lift $\phi$ on a $p$-torsion free ring $S$ (i.e., $p$ is a non-zero divisor in $S$)
one can define the operator $\d:S\ra S$,
$$\d u:=\frac{\phi(u)-u^p}{p},\ \ \ u\in S,$$
which is referred to as the $p$-{\it derivation} attached to $\phi$; following \cite{char,book}
we view $\d$ as an analogue of a derivation.

Similarly if $X$ is a scheme or a $p$-adic formal scheme a 
{\it Frobenius lift} on $X$  is an endomorphism $\phi:X\ra X$ whose reduction mod $p$ 
is the $p$-power Frobenius on the reduction of $X$ mod $p$. 
If $X=Spec\ S$ or $X=Spf\ S$ we usually denote by same letter $\phi$ the Frobenius 
lifts on $S$ and $X$. Some of our rings and schemes to be considered later will be  non-Noetherian; 
typically rings of polynomials in infinitely many variables will have to be considered.
However the only non-Noetherian $p$-adic formal schemes we will encounter will be  affine.

In our paper a special role will be played by the following  ring.  We let $R$ be any complete discrete valuation ring with maximal ideal generated by $p$ and algebraically closed
residue field $\overline{R}:=R/pR$; 
such an $R$ is uniquely determined up to isomorphism by its residue field (it is the Witt ring on that field) and possesses a unique Frobenius 
lift $\phi:R\ra R$.
For any $R$-algebra $S$ and any scheme or $p$-adic formal scheme $X$ over $R$, Frobenius lifts
on $S$ or $X$ will be assumed compatible with the Frobenius lift on $R$. For any ring $S$ and  scheme $X$  we denote by $\widehat{S}$ and $\widehat{X}$ 
the $p$-adic completions of $S$ and $X$ respectively. 
We also define the K\"{a}hler differentials on the affine formal scheme $\widehat{X}$ over a $p$-adically complete ring $S$ by
$$\Omega_{\widehat{X}/S}:=\lim_{\leftarrow} \Omega_{X_n/S_n}$$
where $S_n=S/p^nS$, $X_n=X\otimes S_n$, and $\Omega_{X_n/S_n}$ are the usual K\"{a}hler differentials. If $X$ is smooth over a torsion free $p$-adically complete ring $S$ and $\phi$ is a Frobenius lift on $\widehat{X}$ compatible with a Frobenius lift on $S$ then $\phi$  induces an additive map
\begin{equation}
\label{kedd}
\frac{\phi^*}{p}:\Omega_{\widehat{X}/S}\ra \Omega_{\widehat{X}/S},\ \ \ \omega\mapsto \frac{\phi^*\omega}{p},\end{equation}
which we can think of as a ``normalized" action of $\phi$ on forms, and which,
following \cite{BYM}, we  view  as an analogue of the Lie derivative on forms in classical differential geometry. 

By the way, the operator \ref{kedd} is a lift to characteristic zero of the Cartier operator and plays a key role in the $p$-adic theory of differential forms,  in particular in the theory of the de Rham-Witt complex and in Mochizuki's theory \cite{mochizuki}. 

Alternatively, if one considers an $S$-derivation $\epsilon$ on $\cO_{\widehat{X}}$
then one has a natural additive map 
$$\frac{1}{p}\epsilon\circ \phi:\cO_{\widehat{X}}\ra \cO_{\widehat{X}},\ \ u\mapsto \frac{\epsilon(\phi(u))}{p},$$
hence, for any $\lambda\in S$ one can define the $\lambda$-{\it commutator}
$$[\epsilon,\phi]_{\lambda}:=\frac{1}{p}\epsilon\circ \phi-\lambda\cdot \phi \circ \epsilon:\cO_{\widehat{X}}\ra \cO_{\widehat{X}}.$$
Then $[\epsilon,\phi]_{\lambda}$ is a $\phi$-derivation vanishing on $S$; here, by a $\phi$-{\it derivation} we understand an additive map $D$ satisfying
$$D(uv)=\phi(u)D(v)+\phi(v)D(u).$$
Also, for a $\phi$-derivation $D$ will write $D\equiv 0$ mod $p^n$ if and only if $Du\equiv 0$ mod $p^n$ for all $u$.

\subsection{$\d$-modular functions}
We review here standard terminology from \cite{difmod} and add some new concepts.
It is natural to work over the ring $\bZ_{(p)}$, the localization of $\bZ$ at $(p)=p\bZ$.
We begin by discussing classical modular forms.
Let $z_4,z_6,x,y$ be variables. 
Define the {\it discriminant}  polynomial
$$\Delta=\Delta(z_4,z_6):=4z_4^3+27z_6^2\in \bZ_{(p)}[z_4,z_6];$$
also define the {\it Hasse} polynomial
$$H_p=H=H(z_4,z_6)\in \bZ_{(p)}[z_4,z_6],$$
as the coefficient of $x^{p-1}$ in $(x^3+z_4x+z_6)^{\frac{p-1}{2}}$.
View $\bZ_{(p)}[z_4,z_6]$ as a graded ring with $z_4$ and $z_6$ weighted homogeneous of degree $4$ and $6$ respectively; then $\Delta$ and $H$ are weighted homogeneous of degrees $12$ and $p-1$,
respectively. Let $\Sigma\in \bZ_{(p)}[z_4,z_6]$ be weighted homogeneous, with $\Sigma\not\equiv 0$ mod $p$.
We may then consider the graded rings,
\begin{equation}
 \label{ord}
M:=\bZ_{(p)}[z_4,z_6,\Delta^{-1}],\ \ \ M_{\Sigma}:=\bZ_{(p)}[z_4,z_6,\Delta^{-1}\Sigma^{-1}].\end{equation}
Note that $\widehat{M}\subset \widehat{M_{\Sigma}}$.
Finally recall the {\it $j$-invariant} which is a weighted homogeneous element of degree $0$ of $M$,
$$j:=1728\frac{4z_4^3}{\Delta}\in M.$$

The interpretation of the rings $M,M_H$ (which we will not need) is as follows. 
The scheme $Spec\ M$ is the moduli scheme of pairs consisting of a smooth projective elliptic curve  over a $\bZ_{(p)}$-algebra $S$, equipped with an invertible $1$-form $\omega$ (i.e., a basis for the module of $1$-forms); more precisely any such pair  is isomorphic to a unique pair consisting of the smooth projective model  of  an affine plane curve $y^2=x^3+ax+b$, with $a,b\in S$, equipped with the form $\frac{dx}{y}$.
In particular  the elements of the ring $M$ can be interpreted as {\it modular functions} over $\bZ_{(p)}$. The formal scheme
$Spf\ \widehat{M_H}$ is the {\it ordinary locus} in $Spf\ \widehat{M}$.

Let $S$ be a $\bZ_{(p)}$-algebra.
We will say that $a,b\in S$ is a {\it non-singular pair} if 
$$\Delta(a,b)\in S^{\times}.$$
For $\Sigma\in M$ weighted homogeneous
we will say that $a,b\in S$ is a {\it $\Sigma$-non-singular pair} if it is non-singular and 
$$\Sigma(a,b)\in S^{\times}.$$
A $H$-non-singular pair will simply be called an {\it ordinary pair}.

We now recall a generalization of modular functions introduced in \cite{difmod, book}.
Let us consider,  in addition,  variables $z'_4,z'_6,...,z_4^{(r)}, z_6^{(r)}$. 
Recall from loc.cit. that by 
a {\it $\d$-modular function} of order $r$
one understands an element of the ring $\widehat{M^r}$ where
$$M^r:=M[z'_4,z'_6,....,z_4^{(r)}, z_6^{(r)}].$$
The examples we shall be interested have  $r\leq 1$; but higher order will automatically be brought  into the picture.
Similarly by
a {\it  $\Sigma$-singular $\d$-modular function} of order $r$ one understands an element of the ring $\widehat{M^r_{\Sigma}}$ where
$$M^r_{\Sigma}:=M_{\Sigma}[z'_4,z'_6,....,z_4^{(r)}, z_6^{(r)}].$$
Finally consider the (non-Noetherian!) rings
$$M^{\infty}:=\bigcup_r M^r,\ \ \ M_{\Sigma}^{\infty}:=\bigcup_r M^r_{\Sigma}.$$
There  is a  unique Frobenius lift which we refer to as the {\it universal} Frobenius lift,
\begin{equation}
\label{27mai}
\phi:\widehat{M^{\infty}_{\Sigma}}\ra \widehat{M^{\infty}_{\Sigma}},\end{equation}
such that 
$$\phi(z_4)=z_4^p+pz_4',\ \ \phi(z_6)=z_6^p+pz'_6,\ \ \phi(z'_4)=(z'_4)^p+pz_4'',\ \ \phi(z'_6)=(z'_6)^p+pz''_6,\ \ etc.$$

\begin{remark}

\

1) The term {\it $\Sigma$-singular} is not an adjective: $\d$-modular functions are $\Sigma$-singular $\d$-modular functions and not vice versa.  

2) In \cite{difmod} what we call here ``$\Sigma$-singular" was called ``holomorphic outside $\Sigma$"; changing terminology here makes some of our statements easier to formulate.

3) For $\Sigma=H$ (respectively $\Sigma=1$)
we will say ``ordinary $\d$-modular function", 
(or ``$\d$-modular function", respectively) instead of ``$\Sigma$-singular $\d$-modular function". 

\end{remark}

In view of our applications we will restrict attention, in the discussion below, to $\d$-modular functions of order $1$.

The rings $M^1_{\Sigma}$  carry  a structure  of $\bZ$-graded rings for which the homogeneous elements of degree $k$ are called of {\it weak weight} $k$: the grading  is defined
by letting $z_4,z_6,z'_4,z'_6$ have weak weights $4,6,4p,6p$ respectively.

There is a more refined concept for functions in $\widehat{M^1_{\Sigma}}$ which we refer to as {\it weight} (rather than {\it weak weight}) and which we now recall. This concept is not, strictly speaking, essential for understanding our main result; however some of the questions we would like to raise involve this concept.

First we introduce the following notation. For any $p$-adically complete $p$-torsion free ring $S$ equipped with a Frobenius lift $\phi=\phi^S$,  any  $F=F(z_4,z_6,z'_4,z'_6)\in \widehat{M^1_{\Sigma}}$,  and any $\Sigma$-non-singular pair $a,b \in S$,
we will abusively write
$$F(a,b):=F(a,b,\d a, \d b).$$
 
Now let $\bZ[\phi]$ be the subring generated by $\phi$ in the ring of additive endomorphisms
of our ring $R$ introduced earlier.
For $c\in R^{\times}$ and $w:=\sum a_i\phi^i\in \bZ[\phi]$ we set
$$c^w:=\prod \phi^i(c)^{a_i}\in R^{\times}.$$
An element  $F\in\widehat{M^1_{\Sigma}}$ 
is said to have {\it  weight} $w\in \bZ[\phi]$
if
\begin{equation}
\label{waitfor9}F(c^4a,c^6b)=c^w F(a,b),
\end{equation}
for all $\Sigma$-non-singular pairs $a,b\in R$ and all $c\in R^{\times}$.
An element  $F \in \widehat{M^1_{\Sigma}}$ 
is said to have {\it  weight $w\in \bZ[\phi]$ mod $p^m$}
if
\begin{equation}
\label{waitfor99}F(c^4a,c^6b)\equiv c^w F(a,b)\ \ \ \text{mod}\ \ \ p^m
\end{equation}
for all $\Sigma$-non-singular pairs $a,b\in R$ and all $c\in R^{\times}$. Note some general facts:

1) An element of $\widehat{M^1_{\Sigma}}$  that is congruent mod $p^m$ to an
element of $\widehat{M^1_{\Sigma}}$ 
 of weight $w$ has weight $w$ mod $p^m$. But an element $\widehat{M^1_{\Sigma}}$ 
 of   weight $w$ mod $p^m$ is not a priori  congruent mod $p^m$ to an element of $\widehat{M^1_{\Sigma}}$ 
 of weight $w$.
  
 2) If
 $F\in M_{\Sigma}$ then $F$, viewed as an element of $\widehat{M^1_{\Sigma}}$ has weight $k\in \bZ\subset \bZ[\phi]$  if and only if it is weighted homogeneous of degree $k$. 
 
 3) An element $F\in M^1_{\Sigma}$ has weak weight $k$ if and only if the equality 
\begin{equation}
F(c^4a,c^6b)=c^k F(a,b),
\end{equation}
  holds for all $\Sigma$-non-singular pairs $a,b$ and all roots of unity $c\in R$.
  So if an element $F$ of $M^1_{\Sigma}$ has weight $w=\sum a_i\phi^i$ then, trivially,  $F$ also has weak weight $w(p):=\sum a_ip^i$. Note however that 
   the most interesting $\Sigma$-singular $\d$-modular functions $F\in \widehat{M^1_{\Sigma}}$ in our applications in \cite{difmod,book} do not belong to the ring  $M^1_{\Sigma}$.
   
   4) If an element $F$ of $\widehat{M^1_{\Sigma}}$ has weight $w$ mod $p^m$ then 
   one can show that $F$ is congruent mod $p^m$ to an element of $M^1_{\Sigma}$ that has weak weight $w(p)$. So any element of $\widehat{M^1_{\Sigma}}$ of weight $w$ is a $p$-adic limit of
   elements of $M^1_{\Sigma}$ of weak weight $w(p)$.
   
   \medskip
 
 Finally  the following concept will be involved in our main result.
 A $\Sigma$-singular
$\d$-modular function $F\in M^1_{\Sigma}$ is called {\it quasi linear}
if it is an $M_{\Sigma}$-linear combination of $1, z'_4, z'_6$.
(The terminology is motivated by the classical theory of differential equations where
``quasi linear" means ``linear in the derivatives".)
Let us say that a quasi linear  $F$ as above 
 is  {\it tangential}
if it is  an $M_{\Sigma}$-linear combination of $1, pz'_4, pz'_6$.
If a quasi linear $F$ has weak weight $k$ then it
has the form
\begin{equation}
\label{popo}
F:=\Gamma_k+\Gamma_{k-4p}z'_4+\Gamma_{k-6p}z'_6,\end{equation}
where $\Gamma_k, \Gamma_{k-4p},\Gamma_{k-6p}\in M_{\Sigma}$ are weighted homogeneous of degrees $k, k-4p, k-6p$ respectively. If in addition $F$ is tangential
then $\Gamma_{k-4p}=p\Gamma^*_{k-4p}$, $\Gamma_{k-6p}=p\Gamma^*_{k-6p}$, for some 
$\Gamma^*_{k-4p}, \Gamma^*_{k-6p}\in M_{\Sigma}$.

The link between weak weights and weights mod $p$ or mod $p^2$ is particularly clean in the case of quasi linear functions. Indeed  we have:

\begin{remark}
\label{cryterion}
Let $F$ be  a quasi linear $\Sigma$-singular
$\d$-modular function of weak weight $k$ as in  \ref{popo}. Then the following hold:

1) The function  $F$  has weight   $k+p-\phi$ mod $p$ if and only if
\begin{equation}
\label{ifofu} 4z_4^p\Gamma_{k-4p}+6z_6^p\Gamma_{k-6p}\equiv 0\ \ \ \text{mod}\ \ p.\end{equation}

2) Assume $F$ is tangential. The function  $F$  has weight   $k+p-\phi$ mod $p^2$ if and only if
\begin{equation}
\label{ifof}\Gamma_k+4z_4^p\Gamma^*_{k-4p}+6z_6^p\Gamma^*_{k-6p}\equiv 0\ \ \ \text{mod}\ \ p.\end{equation}
\end{remark}

\begin{example}
Here are some examples of $\Sigma$-singular  quasi linear  $\d$-modular functions that will play a role later:

\medskip

\begin{equation}
\label{archu1}1-p\frac{z'_4}{4z_4^p}\in M^1_{z_4};\end{equation}

\begin{equation}
\label{archu2}1-p\frac{z'_6}{6z_6^p}\in M^1_{z_6};\end{equation}

\begin{equation}
\label{archu3}
1-p\frac{2z_4^{2p}z'_4+9z_6^pz'_6}{2\Delta^p}\in M^1;\end{equation}

\begin{equation}
\label{archu44}
\frac{2z_4^p z'_6-3 z_6^pz'_6}{\Delta^p}\in M^1.\end{equation}

\medskip 

In examples \ref{archu1}, \ref{archu2}, \ref{archu3}, \ref{archu44} we take $\Sigma$ to be
$z_4, z_6, 1, 1$ respectively.
The examples \ref{archu1}, \ref{archu2}, \ref{archu3} are tangential, of weak weight $0$,  and weight $p-\phi$ mod $p^2$. Example \ref{archu44} is not tangential, of weak weight $-2p$, and weight $-p-\phi$ mod $p$. 
By the way
\ref{archu1} and \ref{archu2} can be obtained from \ref{archu3} by setting $z_6=0$ and $z_4=0$ respectively.
\end{example}

\subsection{Elliptic curves} 
We start by reviewing some standard terminology and notation; cf. \cite{silverman}.
In what follows $S$ is, again, a $p$-torsion free $p$-adically complete ring.
For any $a,b\in S$ consider the affine  curve
$$E_{ab}:=Spec\ \frac{S[x,y]}{(y^2-f(x))},\ \ \ f(x):=x^3+ax+b\in S[x].$$
The notation $E_{ab}$ does not specify $S$ which will always be clear from context.
 The condition that the pair $a,b$ be non-singular is equivalent to  $E_{ab}$ being smooth  over $S$ (so an affine elliptic curve); for $S=R$ the condition that the pair $a,b$ be ordinary  is equivalent to $E_{ab}$ being {\it ordinary} (i.e. being with good, ordinary reduction). For $S=R$ and $a,b$ non-singular the condition $a\equiv 0$ mod $p$ is, of course, equivalent to the $j$-invariant satisfying $j\equiv 0$ mod $p$; the condition $b\equiv 0$ mod $p$ is equivalent to $j\equiv 1728$ mod $p$.

For $E_{ab}$ smooth
consider the invertible $1$-form on $E_{ab}$,
$$\omega:=\frac{dx}{y};$$
clearly $\omega$ is a basis  of $\Omega_{\cO(E_{ab})/S}$ and hence of 
$\Omega_{\widehat{E_{ab}}/S}$; $\omega$ extends to the 
smooth projective model ${\mathcal E}_{ab}$ of $E_{ab}$ and 
is  a basis for the global (equivalently, translation invariant) $1$-forms on 
${\mathcal E}_{ab}$.
 Also one can consider the  $S$-derivation on $\cO(E_{ab})$,
$$\epsilon:=y\frac{d}{dx}.$$
This derivation extends to ${\mathcal E}_{ab}$ and is a basis for  the $S$-module of  translation invariant derivations on ${\mathcal E}_{ab}$.

We will be interested in what follows in the open set $E_{ab}^*\subset E_{ab}$ where $y$ is invertible:
$$E_{ab}^*:=Spec\ \frac{S[x,y,y^{-1}]}{(y^2-f(x))}.$$
Invariantly $E_{ab}^*$ can be described as $E_{ab}$  minus its $3$
sections which, on the smooth projective model  ${\mathcal E}_{ab}$
of $E_{ab}$, have order $2$; equivalently, $E^*_{ab}$  is the complement in ${\mathcal E}_{ab}$ of the $2$-torsion sections.
Now $E_{ab}$ carries an involution $[-1]$ induced by the multiplication by $-1$ 
on ${\mathcal E}_{ab}$; on the level of rings the involution is given by $x\mapsto x$, $y\mapsto -y$. 
The involution preserves $E_{ab}^*$. The quotient $E_{ab}/[-1]$ identifies with the affine line $Spec\ S[x]$ and the quotient 
$E^*_{ab}/[-1]$ identifies with the open set of this affine line given by $Spec\ S[x]_f$.
  Moreover
the projection map
$$E_{ab}^*\ra E^*_{ab}/[-1]$$
is \'{e}tale.

\begin{remark}\label{rem1}
Let $E^0_{ab}\subset E_{ab}$ be a Zariski open set.
(Later we will be interested in two cases: $E^0_{ab}=E_{ab}$ and $E^0_{ab}=E^*_{ab}$.)
Consider a Frobenius lift $\phi$ on $\widehat{E_{ab}^0}$, extending a Frobenius lift on $S$,  consider an element $\lambda\in S$, let $n\geq 1$, and consider the following conditions:

\medskip

1) $\phi$ extends to a Frobenius lift 
on the smooth projective model ${\mathcal E}_{ab}$ of $E_{ab}$;

2) $\phi$  satisfies the following congruence in $\Omega_{\widehat{E_{ab}^0}/S}$:
$$\frac{\phi^*}{p}\omega\equiv  \lambda \cdot \omega\ \ \text{mod}\ \ \ 
p^n;$$

3) $\phi$ satisfies the following congruence:
$$[\epsilon,\phi]_{\lambda}\equiv 0 \ \ \ \text{mod}\ \ \  p^n.$$

\medskip

\noindent It is a trivial exercise to check that 
$2\Leftrightarrow 3$. 
Also if 1 holds then 2 holds for some $\lambda$.
By the way, for $S=R$, we 
 also have that condition 1 implies that the smooth projective model of
 $E_{ab}$ is a {\it canonical lift} in the sense of Serre and Tate \cite{Katz};
 and conversely any canonical lift possesses a $\phi$ satisfying condition 1.
 Cf. the Appendix of \cite{messing}. 
\end{remark}

We now introduce the following terminology:

\begin{definition}\label{inv}
Let $a,b\in S$ be a non-singular pair and $\lambda\in S$.
A Frobenius lift $\phi$ on $\widehat{E_{ab}^0}$ is   {\it Lie invariant mod $p^m$} with {\it eigenvalue $\lambda$} if either of the equivalent conditions 2 or 3 in Remark \ref{rem1} holds.
A Frobenius lift $\phi$ on $\widehat{E_{ab}^0}$ is  {\it Lie invariant mod $p^m$} if
it is Lie invariant mod $p^m$, with some eigenvalue.
\end{definition}

As  mentioned in the Introduction  another concept of invariance of Frobenius lifts was introduced in \cite{foundations}, p. 152.

\begin{remark}
Assume, in this remark only, that $S$ is the complex field ${\mathbb C}$ (rather than our $p$-adically complete ring we have been considering so far). Let
 $\d$ be a ${\mathbb C}$-derivation on the ring $\cO(E_{ab}^0)$ (i.e., a {\it vector field} on $E_{ab}^0$)
and let $\text{Lie}_{\d}$ be the induced ${\mathbb C}$-linear endomorphism of $\Omega_{E_{ab}^0/{\mathbb C}}$  (the
{\it Lie derivative} along $\d$). Then the following are equivalent:

\medskip

1) $\d$ extends to a vector field on the smooth projective model ${\mathcal E}_{ab}$ of $E_{ab}$;

2) $\text{Lie}_{\d} \omega=0$;

3) $[\epsilon,\d]=0$;

4) $\d$ is a  ${\mathbb C}$-multiple of $\epsilon$.

\medskip

We would like to see conditions 1, 2, 3 in Remark \ref{rem1} as arithmetic analogues
of the conditions 1, 2, 3 in the present Remark, respectively. Of course condition 4 in the present Remark
has no direct arithmetic analogue.
\end{remark}

Let us close our discussion here by considering elliptic curves over rings of $\d$-modular forms. 
Set ${\mathbb S}:=\widehat{M^{\infty}_{\Sigma}}$,
consider the polynomial 
$$f(z_4,z_6,x):=x^3+z_4x+z_6\in M[x]\subset {\mathbb S}[x],$$
 and consider the elliptic curve
$$
E_{\mathbb S}  :=  Spec\ {\mathbb S}[x,y]/(y^2-f(z_4,z_6,x)),$$
 and its open set
$$E^*_{\mathbb S}:= Spec\ {\mathbb S}[x,y,y^{-1}]/(y^2-f(z_4,z_6,x)).$$
So according to our previously introduced notation $E_{\mathbb S}$ is the curve $E_{z_4z_6}$ over ${\mathbb S}$.
Let $E^0_{\mathbb S}$ be either $E_{\mathbb S}$ or $E^*_{\mathbb S}$.
Assume we are given a  Frobenius lift
$\Phi$ on $\widehat{E^0_{\mathbb S}}$ 
extending the universal Frobenius lift on ${\mathbb S}$ in \ref{27mai} and assume we are given a 
  $\Sigma$-non-singular pair $a,b\in R$. 
  We denote by $E_{ab}$ the elliptic curve over $R$ attached to $a,b$ and we let
  $E_{ab}^0$ be either $E_{ab}$ or $E^*_{ab}$ respectively.
  Then $\Phi$ induces
a Frobenius lift $\phi=\phi_{ab}$ on $\widehat{E^0_{ab}}$: it is  
 the unique morphism, extending the Frobenius lift on $R$, that makes the following diagram commute,
 $$\begin{array}{ccc}
 \widehat{\cO(E^0_{\mathbb S})} &  \stackrel{\Phi}{\longrightarrow} &
\widehat{\cO(E^0_{\mathbb S})}\\
  \downarrow & \ & \downarrow\\
 \widehat{\cO(E^0_{ab})} & \stackrel{\phi_{ab}}{\longrightarrow}
    &  \widehat{\cO(E^0_{ab})}
 \end{array}$$
 where the vertical arrows are  induced by 
 the homomorphism ${\mathbb S}\ra R$ sending 
 $$z_4\mapsto a,\ \ z_6\mapsto b,\ \ z'_4\mapsto \d a,\ \ z'_6\mapsto \d b,\ \ z''_4\mapsto \d^2 a,\ \ z''_6\mapsto \d^2 b,...$$

\section{Main results}

 We will freely use the elliptic curve notation and terminology introduced in the previous section.
First we have the following existence result for   Frobenius lifts that are Lie invariant  mod $p$:

\begin{thm}\label{thm1}
Let ${\mathbb S}:=\widehat{M^{\infty}_H}$. There exists a Frobenius lift  $\Phi$ on $\widehat{E_{\mathbb S}}$, 
   extending the universal Frobenius lift on ${\mathbb S}$ such that $\Phi$ is Lie invariant mod $p$
   with eigenvalue $H^{-1}$.
   \end{thm}
   
   In particular in the situation above we have:
   
   \begin{corollary}\label{tzu}
For any ordinary pair $a,b\in R$  
    the   Frobenius lift $\phi_{ab}$ on $\widehat{E_{ab}}$ induced by
  $\Phi$
 is Lie invariant mod $p$, with eigenvalue $H(a,b)^{-1}$.
\end{corollary}

By the way, conversely, we have:

\begin{proposition}\label{tze}
Assume $a,b\in R$ is a non-singular pair and assume there is a Frobenius lift $\phi_{ab}$ on $\widehat{E_{ab}^*}$ that is Lie invariant mod $p$, with eigenvalue $\lambda\in R$. Then $a,b$ is ordinary and 
$$\lambda\equiv H(a,b)^{-1}\ \ \text{mod}\ \ p.$$
\end{proposition}

Our next Theorem can be viewed as  a lift mod $p^2$  of Theorem \ref{thm1}. 

   \begin{thm}\label{godem5}
   Assume $p\geq 13$,   $p\not\equiv 11$ mod $12$. 
    There exist a triple $\Sigma,\Lambda,\Phi$, where:
    
    i) $\Sigma\in (H)\subset \bZ_{(p)}[z_4,z_6]$ is a weighted homogeneous polynomial,
    $\Sigma\not\equiv 0$ mod $p$;
    
    ii)  $\Lambda \in M^1_{\Sigma}$ is a $\Sigma$-singular $\d$-modular function
     that is quasi linear,  tangential,
   of weak weight $1-p$,
    with $\Lambda\equiv H^{-1}$ mod $p$;
    
    iii)  $\Phi$ is a Frobenius lift 
    on $\widehat{E^*_{\mathbb S}}$ 
   extending the universal Frobenius lift on ${\mathbb S}:=\widehat{M^{\infty}_{\Sigma}}$;
   
   \noindent  such that 
   $\Phi$ is Lie invariant mod $p^2$ with eigenvalue $\Lambda$.\end{thm}

   In particular, in the situation  above, we have:
   
   \begin{corollary}\label{cookcook}
   For any $\Sigma$-non-singular pair $a,b\in R$ 
     the 
  Frobenius lift $\phi_{ab}$ on $\widehat{E_{ab}^*}$ induced by $\Phi$  is Lie invariant mod $p^2$, with eigenvalue
  $\Lambda(a,b)$. \end{corollary}
  
      It is not clear whether one should expect that Theorem \ref{godem5} and Corollary 
      \ref{cookcook} hold with $E^*$ replaced by $E$.
  On the other hand our proof will yield a more precise result.  Indeed, $\Sigma$ and $\Lambda$ in Theorem 
  \ref{godem5}
  can be chosen such that the following Propositions
  \ref{supplement} and \ref{supplement2} hold.

  \begin{proposition}\label{supplement}
  We have $\Sigma=z_6\Psi H$, where $\Psi\in \bZ_{(p)}[z_4,z_6]$ is weighted homogeneous of degree $d$ such that  
  $$\begin{array}{rllllll}
  d&=&p+5 & \ &\text{if} & p\equiv 7  & \text{mod}\ \ 12;\\
 d&=&p+11 & \ & \text{if} & p\equiv 1,5\ \ & \text{mod}\ \ 12;\end{array}$$
 and
 $$\begin{array}{rllllll} \Psi(0,1) & \not\equiv &  0  &  \text{mod}\ \  p &  \text{if} & p\equiv 1& \text{mod}\ \ 3;\\
  \Psi(1,0) & \not\equiv&  0 &  \text{mod}\ \  p &  \text{if} &  p\equiv 1& \text{mod}\ \  4.\end{array}$$ \end{proposition}

 For the next Proposition let us consider the quasi linear ordinary $\d$-modular function of weak weight $1-p$ and weight $1-\phi$ mod $p^2$,
  \begin{equation}
  \label{lambda1}
  \Lambda_1:=\frac{1}{H}\left(1-p\frac{2z_4^{2p}z'_4+9z_6^pz'_6}{2\Delta^p}\right).
  \end{equation}
  
  \begin{proposition}
  \label{supplement2}
 Assume $p\equiv 1$ mod $3$.  There exists $\beta\in \bZ_{(p)}$ such that:
  \begin{equation}
  \label{genev}
  \Lambda \equiv (1-p\beta)\Lambda_1 \equiv  \frac{1}{H}\left(1-p\frac{z'_6}{6z_6^p}-p\beta\right)
  \ \ \ \text{mod}\ \ (p^2,z_4)\ \ \ \text{in}\ \  M^1_{\Sigma}.\end{equation}
  \end{proposition}
  
   Note that if $p\equiv 1$ mod $3$ then
  $z_4$ does not divide $H$ and,
   by Proposition \ref{supplement}, $z_4$ does not divide $\Psi$ either; hence $z_4$ is not invertible in $M^1_{\Sigma}$ and so Proposition \ref{supplement2} is not a tautology.
 On the other hand, by Proposition \ref{supplement}, $z_6$ divides $\Sigma$ in Theorem \ref{godem5} so no pair $a,b$ with $b\equiv 0$ mod $p$  in that Theorem is $\Sigma$-non-singular, so   the Theorem does not apply to such pairs; however, for such pairs, we will prove the following separate result:
  
 \begin{proposition}
 \label{5trays}
  Assume $p\equiv 1$ mod $4$. 
   There exists $\alpha\in \bZ_{(p)}$ such that the following holds: for any  $a,b\in R$ with $a\not\equiv 0$ and $b\equiv 0$ mod $p$
     there exists a  
  Frobenius lift $\phi$ on $\widehat{E_{ab}^*}$ such that $\phi$ is Lie invariant mod $p^2$, with eigenvalue $(1-p\alpha)\Lambda_1(a,b)$.\end{proposition}
  
  It is interesting to compare  Proposition \ref{supplement2} to Proposition \ref{5trays} in conjunction with the fact that:
  \begin{equation}
  \label{vint}
 (1-p\alpha) \Lambda_1\equiv 
  \frac{1}{H}\left(1-p\frac{z'_4}{4z_4^p}-p\alpha\right) \ \ \ \text{mod}\ \ (p^2,z_6)\ \ \ \text{in}\ \ M^1_{z_4H}.\end{equation}
  Indeed in both cases $\Lambda_1$ plays a role although, at this point, it is not clear how to  unify these two cases (even conjecturally). 
   Note by the way that, for $p\equiv 1$ mod $4$, $z_6$ does not divide $H$ hence congruence \ref{vint} is a not a tautology.
   
  \begin{remark}
  The condition  that $a,b$ in Theorem \ref{godem5} be $\Sigma$-non-singular excludes only finitely many values
   for the reduction mod $p$ of the $j$-invariant of $E_{ab}$; hence Theorem \ref{comics} in the Introduction is a consequence of Theorem \ref{godem5} above.\end{remark}

      \begin{remark}\label{seattlee}
       $\Sigma, \Lambda, \Phi$ in Theorem \ref{godem5} are not unique; however note that 
     in our proof $\Sigma, \Lambda, \Phi$   will be canonically constructed, in the sense that they are unique under some natural extra assumptions (which we are not going to make explicit here).
      In addition, the  canonically constructed $\Sigma=z_6\Psi H$    can be explicitly computed, for any given $p$, by a recursive procedure. This computation, in the first two cases ($p=13, 17$) allowed by the Theorem, gives that  $\Psi$ is congruent mod $p$ to a constant times the product $\Delta H$. 
      Now note that, by Proposition \ref{supplement}, we have that both $\Psi$ and $\Delta H$ have  degree $p+11$ if $p\equiv 1,5$ mod $12$. So it seems reasonable 
     to ask if 
      $\Psi$ is congruent mod $p$ to a constant times the product $\Delta H$ for 
      (almost) all $p\equiv 1,5$ mod $12$. 
       Similarly, note that
       for $p\equiv 7$ mod $12$ we have that $z_6$ divides $H$ 
       and that, by Proposition \ref{supplement}, we have that both $z_6\Psi$ and $\Delta H$ have  degree $p+11$. So, again, it seems reasonable to
       ask if 
      $z_6 \Psi$ is congruent mod $p$ to a constant times the product $\Delta H$ for (almost) all $p\equiv 7$ mod $12$. It is not clear if positive answers to these questions should be expected.
      Also it is not clear if our Theorem \ref{godem5} can be extended to the case $p\equiv 11$ mod $12$;
      we will be able to check, however, that Theorem \ref{godem5} holds for $p=11$, with $\Sigma=z_6H$. In particular, 
      in the special cases discussed above, we have:
      \end{remark}

\begin{corollary}\label{godem4}
Assume $p=11,13,17$ and $a,b\in R$ is an ordinary pair.
Then there exists a  
  Frobenius lift $\phi$ on $\widehat{E_{ab}^*}$ such that $\phi$ is Lie invariant mod $p^2$.
\end{corollary}

 We end our discussion of elliptic curves by making a series of remarks on possible links with some remarkable  differential modular forms considered in \cite{difmod, book}:

 \begin{remark}
 \ 
 
  1) We recall from \cite{Barcau, book} that there is 
an   ordinary $\d$-modular function, denoted in loc. cit.  by 
$f^{\partial}\in \widehat{M^1_H}$,
  of  weight $\phi-1$, that plays a key role in the theory of $\d$-modular functions in loc. cit. This form satisfies 
  \begin{equation}
  \label{ferdes}
  f^{\partial}\equiv H\ \ \ \text{mod}\ \ p\end{equation}
   in $\widehat{M^1_H}$. So $f^{\partial}$ is invertible in the ring 
$\widehat{M^1_H}$ and its inverse, denoted in loc. cit. by 
$f_{\partial}\in \widehat{M^1_H}$, 
 has weight $1-\phi$ and satisfies 
 $$f_{\partial}\equiv H^{- 1}\ \ \text{mod}\ \ p.$$
One can ask if there is a  relation  between 
  $f_{\partial}$ and the function
   $\Lambda$ in Theorem \ref{godem5}. By the way, there is a notion of {\it $\d$-Fourier expansion} for $\Sigma$-singular $\d$-modular functions (cf. \cite{difmod}) and the $\d$-Fourier expansion of $f^{\partial}$ was shown to be equal to $1$ (cf. \cite{Barcau, book}).

2) We recall from \cite{difmod, book} that there is a   $\d$-modular function, denoted in loc. cit.  by $f^1\in \widehat{M^1}$, of weight $-1-\phi$, that also plays a key role in the theory of $\d$-modular functions. By  \cite{H}  we have
\begin{equation}
\label{coleg2}f^1\equiv f^1_1+f_0\ \ \ \text{mod}\ \ p,\ \ \ 
f^1_1:= H\cdot \frac{2z_4^pz'_6-3z_6^pz'_4}{\Delta^p},\ \ \ f_0\in M,\end{equation}
where $f_0$ is weighted homogeneous and can be viewed as a  (rather complicated but explicit) ``correction term" for $f^1_1$. 
 Note that $f^1_1$ is quasi linear (non-tangential) 
of weak weight $-1-p$ and  weight $-1-\phi$ mod $p$. One can ask if there is any relation
between $f^1$ and the $\Lambda$ in Theorem \ref{godem5}.

3) Recall from \cite{Barcau} that the forms $f^{\partial}$ and $f^1$ referred to above are related by the formula:
\begin{equation}
\label{coleg3}
f^{\partial}=c\cdot \left(72\cdot \phi(z_6)\frac{\partial}{\partial z'_4}-
16\cdot \phi(z_4)^2\frac{\partial}{\partial z'_6}-p\cdot \phi(P)
\right) f^1,\end{equation}
in $\widehat{M^1_H}$,
where $c\in \bZ_p^{\times}$, and $P\in \widehat{M_H}$ is the Ramanujan function (whose Fourier expansion is the normalized Eisenstein series $E_2$). 
The formulae  \ref{coleg2}, \ref{coleg3}, yield information on $f^{\partial}$ mod $p$ (e.g., they yield \ref{ferdes}) but they  give no information on $f^{\partial}$ mod $p^2$ and hence there seems to be  no direct way, at this point, of comparing  $\Lambda$ in Theorem \ref{godem5} and $f_{\partial}$.

   4) The right hand side of the congruence \ref{genev} is an ordinary $\d$-modular function with  weak weight $1-p$  and weight $1-\phi$ mod $p^2$; one can ask if   $\Lambda$ in Theorem \ref{godem5} can be chosen to have these properties. A strong form of a positive answer to this question would be  that $\Lambda$ in Theorem \ref{godem5} can be taken of the form
   $$\Lambda\equiv \Lambda_0 \Lambda_1\ \ \ \text{mod}\ \ p^2$$
   where $\Lambda_1$ is the function in \ref{lambda1} and $\Lambda_0\in M_H$ 
   is a ``correction factor" that has weight (equivalently, degree) $0$. The situation would then be similar to the one in 2) above. But this picture might be too optimistic.
   \end{remark}

\section{Proofs}

We start with a general discussion on  Frobenius lifts on affine elliptic curves. 
Throughout this section $S$ is a $p$-torsion free $p$-adically complete ring equipped with a Frobenius lift $\phi=\phi^S$. The two cases we are most interested in are $S=R$ and $S={\mathbb S}:=\widehat{M^{\infty}_{\Sigma}}$.

 Let $a,b \in S$ be a non-singular pair. 
 Given a Frobenius lift $\phi$ on $\widehat{\cO(E_{ab}^*)}$ extending $\phi^S$
consider the attached $p$-derivation,
 $$\d:\widehat{\cO(E_{ab}^*)}\ra \widehat{\cO(E_{ab}^*)},\ \ \ 
 u\mapsto \d u=\frac{\phi(u)-u^p}{p}.$$
  The set of  Frobenius lifts $\phi$ on $\widehat{E_{ab}^*}$ extending $\phi^S$ is in bijection with the ring $\widehat{\cO(E_{ab}^*)}$;
the bijection attaches to any 
$\phi$ the 
element $Z:=\d x \in \widehat{\cO(E_{ab}^*)}$. Any element 
 $Z\in \widehat{\cO(E_{ab}^*)}$ comes from some unique $\phi$
  because $S[x]\subset \cO(E_{ab}^*)$ is \'{e}tale.   
  By the way, one easily checks that $\phi$ commutes with the involution $[-1]$
   if and only if $Z\in \widehat{S[x]_f}$.
  
Similarly the derivation $\frac{d}{dx}:S[x]\ra S[x]$ lifts to a unique derivation
$$\frac{d}{dx}:\widehat{\cO(E_{ab}^*)}\ra \widehat{\cO(E_{ab}^*)},\ \ \ u\mapsto \frac{du}{dx};$$
this is, again, because $S[x]\subset \cO(E_{ab}^*)$ is \'{e}tale. 

As a notational convention, we will continue to denote  by $y$  the class of $y$ in $\widehat{\cO(E_{ab}^*)}$.
Also for any element $u$ in a $p$-adically complete ring, with $u\equiv 1$ mod $p$, we denote by $u^{1/2}$  the unique square root in that ring which is $\equiv 1$ mod $p$.

For a Frobenius lift $\phi$ on $\widehat{\cO(E_{ab}^*)}$, extending $\phi^S$,  with $\d x=Z\in \widehat{\cO(E_{ab}^*)}$ we have:
$$\frac{\phi^*}{p}\omega=\frac{d(\phi(x))}{p\phi(y)}=
\frac{d(x^p+pZ)}{p\phi(y)}=\frac{x^{p-1}dx+\frac{dZ}{dx}dx}{\phi(y)}$$
in $$\Omega_{\widehat{\cO(E_{ab}^*)}}=\widehat{\cO(E_{ab}^*)}dx.$$
Hence, for $\lambda\in R$ the condition
\begin{equation}
\label{cond1}
\frac{\phi^*}{p}\omega\equiv \lambda \cdot \omega\ \ \ \text{mod}\ \ \ p^m\end{equation}
is equivalent to the following congruence in $\widehat{\cO(E_{ab}^*)}$:
\begin{equation}
\label{violeta}
\frac{dZ}{dx}+x^{p-1}\equiv \lambda\frac{\phi(y)}{y}\ \ \ \text{mod}\ \ p^m.
\end{equation}
On the other hand set
\begin{equation}
\label{defofK}
K(x):=K_{ab}(x):=\frac{1}{p}(x^{3p}+\phi(a)x^p+\phi(b)-(x^3+ax+b)^p)\in S[x],\end{equation}
 let $z$ be one more variable, 
and consider the cubic polynomial in $z$ with coefficients in $R[x]_f$, 
$$G(x,z):= G_{ab}(x,z)\in S[x]_f[z],$$
given by 
\begin{equation}
\label{cubicc}
G(x,z):= 1+\frac{pK(x)}{f(x)^p}+\frac{p(3x^{2p}+\phi(a))}{f(x)^p}z+\frac{3p^2x^p}{f(x)^p}z^2+\frac{p^3}{f(x)^p}z^3.
\end{equation}
 We then have
$$
\frac{\phi(y)}{y} = y^{p-1}\frac{\phi(y)}{y^p}= f(x)^{\frac{p-1}{2}} \left( \frac{\phi(f(x))}{f(x)^p}\right)^{1/2}= f(x)^{\frac{p-1}{2}} G(x,Z)^{1/2}.$$
Hence congruence \ref{violeta} is equivalent to the congruence
\begin{equation}
\label{weierstrasss}
\frac{dZ}{dx}+x^{p-1}\equiv \lambda f(x)^{\frac{p-1}{2}} G(x,Z)^{1/2}\ \ \ \text{mod}\ \ \ p^m.
\end{equation}
Congruence \ref{weierstrasss} is somewhat reminiscent of the classical differential equation satisfied by the Weierstrass elliptic functions: it implies that a quadratic polynomial $(\frac{dZ}{dx}+x^{p-1})^2$ in $\frac{dZ}{dx}$
is (congruent to)  a cubic polynomial $\lambda^2f(x)^{p-1}G(x,Z)$ in $Z$.
Congruence \ref{weierstrasss} looks quite mysterious in general; we will
be able to solve it, though, in the cases $m=1,2$; in these cases the cubic polynomial
becomes congruent to a linear one.
The case $m=1$ will be dealt with next. The case $m=2$ will then occupy us  for the rest 
of the paper.
\medskip

{\it Proof of 
Theorem \ref{thm1}}.   Take in our discussion above 
$$S=\widehat{M^{\infty}_H},\ \ a=z_4,\ \ b=z_6,\ \ \lambda=H^{-1},\ \ m=1.$$
Then  \ref{weierstrasss} becomes:
\begin{equation}
\label{11}
\frac{dZ}{dx}\equiv \lambda f(x)^{\frac{p-1}{2}}-x^{p-1}\ \ \ \ \text{mod}\ \ \ p.
\end{equation}
The right hand side of \ref{11} 
is a polynomial in $x$ of (non-weighted) degree $\leq \frac{3(p-1)}{2}$ for which the coefficient of $x^{p-1}$ vanishes. 
Such a polynomial is the derivative of a polynomial in $S[x]$ hence the equation \ref{11} has a solution  $Z\in S[x]$. 
Define a    Frobenius lift $\phi$ on $\widehat{E^*_S}$ by setting $\d x=Z$, where $\d$ is attached to $\phi$; then $\phi$ is Lie invariant mod $p$ with eigenvalue $\lambda$. 
This $\phi$ does not a priori extend to $\widehat{E_S}$; in what follows we will modify this $\phi$ so that an extension is possible.

Set $Y(x):=\d f\in S[x]$; a direct computation yields:
\begin{equation}
\label{oortt}
Y(x)\equiv K(x)+(3x^2+a)^pZ(x)\ \ \ \text{mod}\ \ p.
\end{equation}
Let $L$ be an algebraic closure of the fraction field of $\overline{S}:=S/pS$ and denote by $\overline{s},\overline{g}$ the images in $\overline{S},\overline{S}[x]$ of any elements $s\in S$, $g\in S[x]$.  Write
$$\overline{f}(x)=(x-e_1)(x-e_2)(x-e_3)\in L[x]$$ with $e_1,e_2,e_3\in L$.
For any $\mu=(\mu_0,\mu_1,\mu_2)\in S^3$  set 
$$Z^{\mu}(x):=Z(x)+\mu_0+\mu_1x^p+\mu_2x^{2p}\in S[x],$$
$$Y^{\mu}(x):=K(x)+(3x^2+a)^pZ^{\mu}(x)\in S[x].$$
We claim that there exist (unique) elements $\overline{\mu_0},\overline{\mu_1},\overline{\mu_2}\in \overline{S}$ (where $\mu_i\in S$)
such that
the following holds in $L$:
\begin{equation}
\label{rufe}
\overline{Y^{\mu}}(e_i)= 0,\ \ \ \text{for}\ \ \ i=1,2,3.\end{equation}
Indeed \ref{rufe} can be viewed as a linear system with coefficients in $L$ and unknowns $\overline{\mu_0},\overline{\mu_1},\overline{\mu_2}\in L$;
the square of the determinant of this system is, up to a sign,  a power of the discriminant of $\overline{f}$  (cf. \cite{lang}, p. 204) and Cramer's rule plus the fundamental theorem of symmetric polynomials easily imply  that the solution to the system has all its components in $\overline{S}$; this proves our claim. Now \ref{rufe} implies that
 $\overline{Y^{\mu}}$ is divisible by $\overline{f}$ 
in $L[x]$ and hence in $\overline{S}[x]$.
Let $\phi^{\mu}$ be the Frobenius lift on $\widehat{E^*}$ with attached $p$-derivation $\delta^{\mu}$ defined by $\d^{\mu} x=Z^{\mu}$. Since 
$$\frac{dZ^{\mu}}{dx}\equiv \frac{dZ}{dx}\ \ \text{mod}\ \ p\ \ \ \text{in}\ \ S[x]$$
it follows that $\phi^{\mu}$ is Lie invariant mod $p$ with eigenvalue $H^{-1}$.
Also $$Y^{\mu}=\d^{\mu}f=\d^{\mu} (y^2)\equiv 2y^p\d^{\mu}y\ \ \ \text{mod}\ \ p,$$
hence
\begin{equation}
\label{maduc}
\d^{\mu}y\equiv \frac{Y^{\mu}}{2y^p}\equiv \frac{Y^{\mu}}{2f^{\frac{p+1}{2}}}y\ \ \text{mod}\ \ p.\end{equation}
 On the other hand a direct computation yields:
$$\begin{array}{rcl}
\frac{d Y^{\mu}}{dx} & \equiv  & \frac{dK}{dx}+(3x^{2}+a)^p \frac{dZ}{dx}\\
\ & \ & \ \\
\ & \equiv & 3x^{3p-1}+a^px^{p-1}-f^{p-1}(3x^2+a)+(3x^{2}+a)^p(\lambda f^{\frac{p-1}{2}}-x^{p-1})\\
\ & \ & \ \\
\ & \equiv & -f^{p-1}(3x^2+a)+(3x^2+a)^p\lambda f^{\frac{p-1}{2}}\ \ \ \text{mod}\ \ p.\end{array}
$$
So $\frac{d\overline{Y^{\mu}}}{dx}$ is divisible by $\overline{f}^{\frac{p-1}{2}}$ in $\overline{S}[x]$. Since $\overline{Y^{\mu}}$  is also divisible by $\overline{f}$ it follows that $\overline{Y^{\mu}}$ is divisible by $\overline{f}^{\frac{p+1}{2}}$ in $\overline{S}[x]$.
By \ref{maduc}, $\d^{\mu}y$ is congruent mod $p$ to an element of $\cO(E_{ab})$.
It follows that $\phi^{\mu}$ induces an endomorphism $\phi^{\mu}_2$ of $E_{ab}\otimes S/p^2S$.
By smoothness of $E_S$ over $S$, $\phi^{\mu}_2$ lifts to an endomorphism $\Phi$ of $\widehat{E_{ab}}$
which will be automatically Lie invariant mod $p^2$ with eigenvalue $\lambda$.
\qed

\bigskip

{\it Proof of Proposition \ref{tze}.}
Assume the hypothesis of the Proposition. Then \ref{weierstrasss} holds for $m=1$, i.e.,
\ref{11} holds.
Now $$\widehat{\cO(E^*_{ab})}=\widehat{R[x]_f}\oplus \widehat{R[x]_f}\cdot y,$$
so one can write
$$Z\equiv Z_0+Z_1y\ \ \text{mod}\ \ p$$
for some $Z_0,Z_1\in R[x]_f$.
We have
$$\frac{dZ}{dx}=\frac{dZ_0}{dx}+\frac{dZ_1}{dx}y+Z_1\frac{dy}{dx}=\frac{dZ_0}{dx}+\left(\frac{dZ_1}{dx}+Z_1\frac{1}{2f}\frac{df}{dx}\right)y.
$$
Hence denoting by an overline the class mod $p$ we have that \ref{11} implies
\begin{equation}
\label{oohh}
\frac{d\overline{Z_0}}{dx}=\overline{\lambda} f(x)^{\frac{p-1}{2}}-x^{p-1}
\end{equation}
in $\widehat{\cO(E^*_{ab})}/(p)$, hence in 
$$R[x]_f/(p)\subset \overline{R}((x)):=\overline{R}[[x]][x^{-1}].$$
Now the left  hand side of \ref{oohh} has no monomial in $x^{p-1}$ and the Proposition follows.
\qed

\bigskip

In  what follows we  continue  our general discussion made before the proofs of Theorem \ref{thm1}
and Proposition \ref{tze}.

By  \ref{weierstrasss},
 the condition that a Frobenius lift $\phi$ on $\widehat{E_{ab}^*}$, extending $\phi^S$, with $\d x=:Z$, and an element $\lambda\in S$ satisfy
\begin{equation}
\label{cond2}
\frac{\phi^*}{p}\omega\equiv \lambda \cdot \omega\ \ \ \text{mod}\ \ \ p^2\end{equation}
is equivalent to the following congruence in $\widehat{\cO(E^*_{ab})}$:
\begin{equation}
\label{www}
\frac{dZ}{dx} \equiv  \lambda f(x)^{\frac{p-1}{2}}- x^{p-1}+
\frac{p\lambda}{2f(x)^{\frac{p+1}{2}}}((3x^{2p}+a^p)Z+K(x))
\ \ \ \ \text{mod}\ \ \ p^2.
\end{equation}  
In what follows, for given $a,b\in S$ with properties to be determined,  we seek a rational function 
$Z=Z(x)\in S[x]_f$ and a $\lambda \in S$ satisfying \ref{www}; for such a $Z(x)$ the formula
 \begin{equation}
 \label{friday}
 \phi(x):=x^p+pZ(x)
 \end{equation}
 defines  a   Frobenius lift $\phi$ on $\widehat{E_{ab}^*}$ 
that is Lie invariant mod $p^2$, with eigenvalue $\lambda$. 

Start with elements   $a,b,\lambda,\theta\in S$ satisfying
\begin{equation}
\label{defth}
\begin{array}{rcl}
\Delta(a,b) & \in & S^{\times},\\
\ & \ & \ \\
 \lambda H(a,b)-1 & \equiv & p\cdot \theta\ \ \ \text{mod}\ \ p^2.\end{array}\end{equation}
In particular the pair $a,b$ is ordinary.
We have:
$$\begin{array}{rcl}
K(x) & \equiv & K_0(x)+\d b+\d a \cdot x^p\ \ \ \text{mod}\ \ p,\ \ \text{where}\\
\ & \ & \ \\
K_0(x) &:=&\frac{1}{p}(x^{3p}+a^px^p+b^p-(x^3+ax+b)^p).
\end{array}$$
So the coefficients of $K_0(x)$ are given by universal polynomials
in $\bZ_{(p)}[z_4,z_6]$ evaluated at $a,b$. (Here and later by ``universal" we mean
``depending on $p$ only".)
For $a,b,\lambda,\theta$ satisfying  \ref{defth},
let 
$$\lambda_0:=H(a,b)^{-1}\in S;$$
so we have
 $$\lambda\equiv \lambda_0(1+p\theta)\ \ \text{mod}\ \ p^2. $$
Let $W(x)\in S[x]$ be the unique polynomial satisfying
 \begin{equation}
 \label{grab}
\frac{dW}{dx}(x)=\lambda f(x)^{\frac{p-1}{2}}-x^{p-1},\ \ \ W(0)=0.\end{equation}
 This $W(x)$ exists because  the right hand side of \ref{grab}
is a sum of monomials of (non-weighted) degree $\not\equiv -1$ mod $p$ plus the monomial
$$(\lambda H(a,b)-1)x^{p-1}=p\theta \cdot x^{p-1}.$$
 Similarly let $W_0(x)\in S[x]$ be the unique polynomial
satisfying
\begin{equation}
 \label{grabber}
 \frac{dW_0}{dx}(x)=\lambda_0 f(x)^{\frac{p-1}{2}}-x^{p-1},\ \ \ W_0(0)=0.\end{equation}
The coefficients of $W_0$ are given by some universal
functions in $M_H$ evaluated at $a,b$, and the coefficient of $x^{p-1}$ in $W_0(x)$ vanishes.
 Then we have
 $$W(x)\equiv  W_0(x)+\theta \cdot x^p\ \ \ \text{mod}\ \ p.$$
 Next we will seek a solution $Z(x)\in S[x]_f$ to \ref{www} of the form
 \begin{equation}
 \label{gogol}
 Z(x)=W(x)+V(x^p)+p\frac{U(x)}{f(x)^p},\end{equation}
 for some polynomials $V(x), U(x)\in S[x]$.
 In terms of the unknown polynomials $V,U$ congruence \ref{www} becomes:

\begin{equation}
\label{ppo}
\frac{dU}{dx}(x)\equiv A(x)+B(x)+C(x)+D(x)+E(x)+F(x)\ \ \ \text{mod}\ \ \ p,\end{equation}
where:

\medskip

 \begin{equation}
 \label{ABCD}
 \begin{array}{rcl}
 A(x) & := & -x^{p-1}f(x)^p\frac{dV}{dx}(x^p),\\
 \ & \ & \ \\
 B(x) & := & \frac{\lambda_0}{2}f(x)^{\frac{p-1}{2}}(3x^{2p}+a^p)V(x^p),\\
 \ & \  & \ \\
 C(x) & := & \frac{\lambda_0}{2}f(x)^{\frac{p-1}{2}}(3x^{2p}+a^p) x^p\cdot \theta,\\
 \ & \ & \ \\
 D(x) & := & \frac{\lambda_0}{2}f(x)^{\frac{p-1}{2}}(K_0(x)+(3x^{2p}+a^p)W_0(x)),\\
 \ & \ & \ \\
 E(x) & := & \frac{\lambda_0}{2}f(x)^{\frac{p-1}{2}}\cdot \d b,\\
 \ & \ & \ \\
 F(x) & := & \frac{\lambda_0}{2}f(x)^{\frac{p-1}{2}}x^p\cdot \d a.\end{array}
 \end{equation}
 
 \medskip
 
 So in order to solve \ref{www} it is enough to find a polynomial $V(x)\in S[x]$ such that
  the
 right hand side of \ref{ppo} has no terms in $x^{sp-1}$ for $s\geq 1$; for if this is the case
 the right hand side of \ref{ppo} is congruent mod $p$ to the derivative of some polynomial $U(x)\in S[x]$
 and we are done.
 
 Write 
 \begin{equation}
 \label{DV}
 V(x)=\sum_{j\geq 0} v_jx^j\end{equation}
 with $v_j$ unknown (and almost all $0$).
 Moreover let $a_s, b_s, c_s, d_s, e_s,f_s$ be the coefficients of $x^{sp-1}$ in $A, B, C, D, E, F$ respectively. Then we have congruences mod $p$:
 
 \medskip
 
 \begin{equation}
 \label{abcd}
 \begin{array}{rcl}
  a_s & \equiv  & -sb^pv_s-(s-1)a^p v_{s-1}-
  (s-3)v_{s-3},\\
 \ & \ & \ \\
 b_s & \equiv  & \frac{1}{2}a^p v_{s-1}+\frac{3}{2}v_{s-3},\\
  \ & \ & \ \\
  c_2  & \equiv   & \frac{1}{2}a^p \theta,\ \ \ c_4\equiv \frac{3}{2}\theta,\ \ \ \text{and}\ \ \ c_s\equiv 0\ \ \text{for}\ \ s\neq 2,4,\\
  \ & \ & \ \\
  e_1 & \equiv  & \frac{\d b}{2}\ \ \ \text{and}\ \ \ e_s\equiv 0\ \ \ \text{for}\ \ s\neq 1,\\
   \ & \ & \ \\
  f_2 & \equiv  & \frac{\d a}{2}\ \ \ \text{and}\ \ \ f_s\equiv 0\ \ \ \text{for}\ \ s\neq 2,\\
   \ & \ & \ \\
   d_s & \equiv & D_s(a,b)\ \ \text{for}\ \ 1\leq s\leq 4\ \ \text{and}\ \ d_s\equiv 0\ \ \text{for}\ \ d\geq 5,
 \end{array}
 \end{equation}
 
 \medskip
 
 \noindent
 for $D_1,D_2,D_3,D_4\in M_H$  some universal  weighted homogeneous elements of degree $6p, 4p, 2p, 0$, respectively. The vanishing of $d_s$ for $s\geq 5$ holds because 
 the (non-weighted) degree of $D(x)$ in $x$ 
 is at most
  $5p-2$.  To check that $D_s$ are weighted homogeneous let us consider  the ring $M_H[x]$ as graded with $z_4,z_6,x$ weighted homogeneous of degree $4,6,2$. Then 
  $$\lambda_0,\ \ f(x),\ \ K_0,\ \ 3x^{2p}+a^p,\ \ W_0(x)$$
   come from universal weighted homogeneous elements  of $M_H[x]$ of degree 
   $$1-p,\ \ 6,\ \ 6p,\ \ 4p,\ \ 2p,$$ respectively. So $D(x)$, and hence $d_sx^{sp-1}$,  come from universal weighted homogeneous elements of degree $8p-2$, hence $d_s$ come from  universal weighted homogeneous elements of degree $$8p-2-2(sp-1)=(8-2s)p.$$

Now the condition that the reduction mod $p$ of the right hand side of \ref{ppo} have no terms in $x^{sp-1}$ for $s\geq 1$
can be written as:
 $$a_s+b_s+c_s+d_s+e_s+f_s\equiv 0\ \ \ \text{mod}\ \ p,\ \ \ s\geq 1,$$
 or, equivalently,
 \begin{equation}
 \label{louvre}
 \begin{array}{rcl}
 sb^pv_s & \equiv  & \left(\frac{3}{2}-s\right)a^p v_{s-1}+
 \left(\frac{9}{2}-s\right)v_{s-3}\\
 \ & \ & \ \\
 \ & \ & +c_s+d_s+e_s+f_s\ \ \ \text{mod}\ \ \ p,\ \ \ s\geq 1.
 \end{array}
 \end{equation}
 
 The following statement summarizes our discussion:
 
  \begin{proposition}
 \label{operastabbing}
 Let
 $a,b,\lambda,\theta$ satisfy \ref{defth} 
 and let 
 $c_s,d_s,e_s,f_s$ 
 be the coefficients of $x^{sp-1}$ in the polynomials $C,D,E,F$ defined in \ref{ABCD}.
 Let 
 ${\mathfrak S}$ be the system of congruences $\ref{louvre}$, for $s\geq 1$, 
 viewed as a system in the unknowns $v_0,v_1,v_2,...$, where $v_{-1},v_{-2}$ are taken to be $0$. 
 Assume $v_0,v_1,v_2,..\in S$ is a solution to the system ${\mathfrak S}$ almost all of whose components are $0$. 
 Define $W$ by the formula \ref{grab}, define
 $V$ by the formula \ref{DV}, define $U$ by the formula \ref{ppo}, define $Z$ by the formula
 \ref{gogol}, and define 
    the Frobenius lift $\phi$ on $\widehat{E^*_{ab}}$ by the formula \ref{friday}. 
    Then $\phi$  is Lie invariant mod $p^2$ with eigenvalue $\lambda$.
\end{proposition}
    
   In what follows we will analyze the system ${\mathfrak S}$.
   For any integer $T$
 let us denote by ${\mathfrak S}_{T}$ the system obtained from ${\mathfrak S}$ by just retaining the congruences \ref{louvre} for $1\leq s\leq T$. Note then that given $a,b,\theta, v_0\in R$ with $b\not\equiv 0$ mod $p$, and given $1\leq T\leq p-1$ there is always a solution 
 to ${\mathfrak S}_{T}$ of the form $v_0,v_1,...,v_{T}\in S$ and this  solution is unique up to congruence mod $p$. We have the following:
 
 \begin{proposition}
 \label{stabcr}
  Assume  $a,b,\theta, v_0\in S$ with $a,b$ ordinary and $b\in S^{\times}$. Assume $p\geq 11$ (hence $\frac{p+7}{2}\leq p-1$)
  and let  $$v_0,v_1,...,v_{\frac{p+7}{2}}\in S$$ 
 be the unique mod $p$ solution 
 to ${\mathfrak S}_{\frac{p+7}{2}}$ attached to $a,b,\theta, v_0$. Assume
 \begin{equation}
 \label{planeplan}
 v_{\frac{p+5}{2}}\equiv  v_{\frac{p+7}{2}}\equiv 0\ \ \ \text{mod}\ \ p.\end{equation}
 Then $$v_0,v_1,...,v_{\frac{p+3}{2}},0,0,...\in S$$ is a solution to ${\mathfrak S}$.
 \end{proposition}
 
 {\it Proof}.
 Congruence \ref{louvre} for $s=\frac{p+9}{2}$ reads
 $$\frac{p+9}{2}b^p v_{\frac{p+9}{2}}\equiv \left(\frac{3}{2}-\frac{p+9}{2}\right)a^pv_{\frac{p+7}{2}}+ \left(\frac{9}{2}-\frac{p+9}{2}\right)v_{\frac{p+3}{2}}\ \ \ \text{mod}\ \ p;$$
this congruence  is satisfied if we set $v_{\frac{p+9}{2}}=0$   because, by \ref{planeplan}, 
  $v_{\frac{p+7}{2}}\equiv 0$ mod $p$ and, on the other hand, 
  $$ \frac{9}{2}-\frac{p+9}{2}=-\frac{p}{2}\equiv 0\ \ \text{mod}\ \ p\ \ \text{in}\ \ \bZ_{(p)}.$$
  Next, the congruence \ref{louvre}  for  $s=\frac{p+11}{2}$ 
  reads
  $$\frac{p+11}{2}b^p v_{\frac{p+11}{2}}\equiv \left(\frac{3}{2}-\frac{p+11}{2}\right)a^pv_{\frac{p+9}{2}}+ \left(\frac{9}{2}-\frac{p+11}{2}\right)v_{\frac{p+5}{2}}\ \ \ \text{mod}\ \ p;$$
this congruence is again satisfied 
  by further setting $v_{\frac{p+11}{2}}=0$ because, by \ref{planeplan}, $v_{\frac{p+5}{2}}\equiv 0$ mod $p$. By induction one gets that \ref{louvre} is satisfied for all $s\geq \frac{p+11}{2}$ by setting $v_s=0$ for all such $s$. \qed
 
 \bigskip

 We would like to explore consequences of the above Proposition. Assume we are given
 an ordinary pair $a,b\in S$ with $b\in S^{\times}$. Then for any $v_0,\theta\in S$
 one can   express the components $v_n$ of the solution 
 $v_0,v_1,...,v_{\frac{p+7}{2}}$ to ${\mathfrak S}_{\frac{p+7}{2}}$ in the form
 $$v_n=\alpha_nv_0+\beta_n\theta+\mu_n \d a +\nu_n \d b +\eta_n$$
 where $\alpha_n,\beta_n, \mu_n, \nu_n, \eta_n$ are given by universal functions in $M_{z_6 H}$ evaluated at $a,b$.  Explicitly, the first values of $\alpha_n,\beta_n$ are given by 
    
    \medskip
    
    $$
    \begin{array}{ll}
    \alpha_1=\frac{a^p}{2b^p}, & \beta_1=0,\\
    \ & \ \\
    \alpha_2=-\frac{a^{2p}}{8b^{2p}}, & \beta_2=\frac{a^p}{4b^p},\\
     \ & \ \\
    \alpha_3= \frac{a^{3p}}{16b^{3p}}+\frac{1}{2b^p}, & \beta_3=-\frac{a^{2p}}{8b^{2p}},\\
     \ & \ \\
\alpha_4=-\frac{5a^{4p}}{128b^{4p}}-\frac{a^p}{4b^{2p}}, & \beta_4=\frac{5a^{3p}}{64b^{3p}}+\frac{3}{8b^p},
    \end{array}
    $$
    
    \medskip
    
   \noindent  while, for $5\leq n\leq p-1$, we have
   
   \medskip
   
    $$
    \begin{array}{rcl}
    \alpha_n & = & -\frac{2n-3}{2n}\frac{a^p}{b^p}\alpha_{n-1}-\frac{2n-9}{2n}\frac{1}{b^p}\alpha_{n-3},\\
    \ & \ & \ \\
     \beta_n & = & -\frac{2n-3}{2n}\frac{a^p}{b^p}\beta_{n-1}-\frac{2n-9}{2n}\frac{1}{b^p}\beta_{n-3}.
    \end{array}
    $$
    
    \medskip

    \noindent For $n\leq p-2$ consider the determinants:
    $$\psi_n:=\text{det} \left(
    \begin{array}{cc}
    \alpha_n & \beta_n\\
    \alpha_{n+1} & \beta_{n+1}\end{array}
    \right).$$
    The first values of these determinants are:
   $$
   \begin{array}{rcl}
   \psi_1 & = & \frac{a^{2p}}{2^3b^{2p}},\\
   \ & \ & \ \\
   \psi_2 & = & -\frac{a^p}{2^3b^{2p}},\\
   \ & \ & \ \\
   \psi_3 & = & \frac{a^{3p}}{2^5b^{4p}}+\frac{3}{2^4b^{2p}}.
   \end{array}
   $$
  Also the above recurrence relations for $\alpha_n,\beta_n$ imply the following  relations:
   \begin{equation}
   \label{tututu}
   \begin{array}{rcl}
   \psi_n & = & \frac{2n-7}{2n+2}\frac{1}{b^p}(\alpha_{n-2}\beta_n-\alpha_n\beta_{n-2}),\ \ 4\leq n\leq p-2,\\
   \ & \ & \ \\
   \psi_n &= &-\frac{2n-7}{2n+2}\frac{2n-3}{2n}\frac{a^p}{b^{2p}}\psi_{n-2}+
   \frac{2n-7}{2n+2}\frac{2n-9}{2n}\frac{1}{b^{2p}}\psi_{n-3},\ \ \ 5\leq n\leq p-2.\end{array}
   \end{equation}
   By induction we get:

   \begin{proposition} \label{cal} For all $n$ with 
   $1\leq n\leq \frac{p-3}{2}$ we have
   $$\psi_{2n}=\frac{\Psi_{2n}(a^p,b^p)}{b^{2np}},\ \ \psi_{2n-1}=\frac{\Psi_{2n-1}(a^p,b^p)}{b^{2np}},$$
   for some universal 
    weighted homogeneous polynomials $\Psi_{2n},\Psi_{2n-1}\in \bZ_{(p)}[z_4,z_6]$ of degrees $4n$ and $4n+4$, respectively. \end{proposition}

   We now prove:

   \begin{thm}\label{godem1}
   Let $p\geq 11$, set $\Psi=\Psi_{\frac{p+5}{2}}$, assume
   $\Psi\not\equiv 0$ mod $p$, and let $\Sigma:=z_6\Psi H$. 
   There exists a $\Sigma$-singular $\d$-modular function, $\Lambda \in M^1_{\Sigma}$,
   that is quasi linear, tangential, of weak weight $1-p$, with $\Lambda\equiv H^{-1}$ mod $p$, 
   and there exists a Frobenius lift $\Phi$ on $\widehat{E^*_{\mathbb S}}$ 
   extending the universal Frobenius lift on ${\mathbb S}:=\widehat{M^{\infty}_{\Sigma}}$, such that 
   $\Phi$ is Lie invariant mod $p^2$ with eigenvalue $\Lambda$.   \end{thm}

   \medskip
   
   Note the hypothesis $\Psi_{\frac{p+5}{2}}\not\equiv 0$ mod $p$;
   conditions when this is satisfied will be provided later (cf. Proposition \ref{golem}).
   
   \medskip

   {\it Proof}.
   Consider two variables $v,t$ and consider the ring
   $${\mathbb M}:=M^1_{\Sigma}[v,t].$$
   Consider the system of linear equations ${\mathfrak S}^{\mathbb M}$ in $v_0,v_1,v_2,......$
   with coefficients in ${\mathbb M}$ obtained from 
   the congruences \ref{louvre} by taking   $a, b,\d a, \d b, \theta$
   to be
   $z_4, z_6, z'_4, z'_6, t$.
   Also let ${\mathfrak S}^{\mathbb M}_{\frac{p+7}{2}}$ be its truncation, as usual.
   The  system ${\mathfrak S}^{\mathbb M}_{\frac{p+7}{2}}$   has a unique (mod $p$) solution, 
   $$v_0, v_1, ... , v_{\frac{p+7}{2}}\in {\mathbb M}$$
   with $v_0=v$. Now consider the graded ring structure on ${\mathbb M}$ defined by letting
   the variables
   $$z_4, z_6, z'_4, z'_6, v, t$$
   have weights
   $$4, 6, 4p, 6p, 2p, 0,$$
   respectively.
   It is trivial to see by induction that $v_n$ is weighted homogeneous of degree
    $2p-2np$ for all $0\leq n\leq \frac{p+7}{2}$. So the ideal $(v_{\frac{p+5}{2}},v_{\frac{p+7}{2}})\subset {\mathbb M}$
    is graded and hence we have a morphism
    of graded rings
    $$M^1_{\Sigma}\ra {\mathbb M}/(v_{\frac{p+5}{2}},v_{\frac{p+7}{2}}).$$
    This morphism is an isomorphism because $\Psi_{\frac{p+5}{2}}$ is invertible in $M^1_{\Sigma}$. Denote by $\tilde{m}\in M^1_{\Sigma}$ the element that maps to the class of an element $m\in {\mathbb M}$.
    So $\Theta:=\tilde{t}\in M^1_{\Sigma}$ has weak weight $0$. This $\Theta$ is trivially seen to be an $M_{\Sigma}$-linear combination of $1,z'_4,z'_6$ so 
    $$\Lambda:=\frac{1}{H}(1+p\Theta)\in M^1_{\Sigma}$$
     is a tangential quasi linear function of weak weight $1-p$. 
     Let ${\mathfrak S}^{\mathbb S}$ 
     be the system with coefficients in $M^1_{\Sigma}\subset {\mathbb S}$
     obtained from ${\mathfrak S}^{\mathbb M}$ by applying the map $m\mapsto \tilde{m}$ to the coefficients
     and let ${\mathfrak S}^{\mathbb S}_{\frac{p+7}{2}}$ be its corresponding truncation. Then 
      \begin{equation}
    \label{zit}
    \tilde{v}_0,\tilde{v}_1,\tilde{v}_2,...,\tilde{v}_{\frac{p+3}{2}},0,0\in M^1_{\Sigma}\end{equation}
    is a solution to
     ${\mathfrak S}^{\mathbb S}_{\frac{p+7}{2}}$. By  Proposition \ref{stabcr} 
      \begin{equation}
    \label{zitt}
    \tilde{v}_0,\tilde{v}_1,\tilde{v}_2,...,\tilde{v}_{\frac{p+3}{2}},0,0,0,...\in M^1_{\Sigma}\end{equation}
    is a solution to
     ${\mathfrak S}^{\mathbb S}$. By Proposition \ref{operastabbing}
     one can attach to \ref{zitt} a
     Frobenius lift $\Phi$ on $\widehat{E^*_{\mathbb S}}$ that is Lie invariant mod $p^2$ with eigenvalue $\Lambda$ and we are done.
       \qed

  \begin{remark}
  What we intuitively did in the above proof was to  view the congruences \ref{planeplan} as  a linear system 
  ${\mathfrak R}$ in the ``unknowns" $v_0$ and $\theta$ with coefficients in $M^1_{\Sigma}$. We then solved the system 
  for $v_0,\theta$ and considered the system ${\mathfrak S}^{\mathbb S}_{\frac{p+7}{2}}$ 
  with these values of $v_0,\theta$. Finally we used a solution to ${\mathfrak S}^{\mathbb S}_{\frac{p+7}{2}}$
  to construct $\Phi$. The system ${\mathfrak R}$ will play again a role later. \end{remark}
        
   \medskip
   
   \begin{remark}
  Here are some further  computations using our recurrence relations \ref{tututu}:
  
  \medskip
  
   $$
   \begin{array}{rcl}
   \psi_4 & = & \frac{a^{2p}}{2^7 \times  5 b^{4p}}\\
   \ & \ & \ \\
   \psi_5 & = & -\frac{7a^{4p}}{2^8 \times 5 b^{6p}}-\frac{23a^p}{2^7 \times  5 b^{4p}}\\
   \ & \ & \ \\
   \psi_6 & = & \frac{17a^{3p}}{2^{10} \times 7 b^{6p}}+\frac{15}{2^7 \times  7  b^{4p}}\\
   \ & \ & \ \\
   \psi_7 & = & \frac{7\times 11a^{5p}}{2^{13}\times 5b^{8p}}+\frac{3\times 43a^{2p}}{2^{11}\times 5b^{6p}}\\
   \ & \ & \ \\
   \psi_8 & = & -\frac{2477 a^{4p}}{2^{15}\times 5 \times 7 b^{8p}}-\frac{2102a^p}{2^{12}\times 5 \times 7 b^{6p}}\\
    \ & \ & \ \\
    \psi_9 & = & -\frac{7\times 11^2 a^{6p}}{2^{16}\times 3\times 5b^{10p}}
    -\frac{2937a^{3p}}{2^{14}\times 5\times 7b^{8p}}+\frac{3\times 11}{2^{10}\times 7b^{6p}},\ \ \text{etc}.
   \end{array}
   $$
   
   \bigskip
   
   \noindent In particular, for $p=11$, we have
   
   \medskip
   
   $$\begin{array}{rcll}
   \Psi_{\frac{11+5}{2}}=\Psi_8 &\equiv &
   z_4(z_4^3+4z_6^2) & \text{mod}\ \ 11,\\
   \ & \ & \ & \ \\
   \Delta  & \equiv  & 4(z_4^3+4z_6^2) & \text{mod}\ \ 11,\\
   \ & \ & \ & \ \\
   H=H_{11} & \equiv & -2 z_4z_6 & \text{mod}\ \ 11,\end{array}$$
   hence
   \begin{equation}
   \label{bb11}
   z_6\Psi_8  \equiv   3\Delta H_{11} \ \ \ \text{mod}\ \ 11.\end{equation}
   
   \medskip
   
     \noindent Also, for $p=13$, we have:
     
     \medskip 
     
    $$\begin{array}{rcll}
   \Psi_{\frac{13+5}{2}}=\Psi_9 & \equiv & 
   4(z_4^6+z_4^3z_6^2+z_6^4)& \text{mod}\ \  13\\
    \ & \ & \ &\ \\
   \  & \equiv
   &  4(z_4^3-9z_6^2)(z_4^3-3z_6^2)&  \text{mod}\ \ 13,\\
   \ & \ & \ &\ \\
   \Delta & \equiv  & 4(z_4^3-3z_6^2)& \text{mod}\ \ 13,\\
    \ & \ & \ &\ \\
   H=H_{13} & \equiv  & 7(z_4^3-9z_6^2)&  \text{mod}\ \ 13,\end{array}$$
    hence
    \begin{equation}
    \label{bb13}
    \Psi_{9}  \equiv  2\Delta H_{13}\ \ \   \text{mod}\ \ 13.\end{equation}
   
   \medskip

  \noindent  Similarly, for $p=17$, we have:
  
     \medskip
   
   $$\begin{array}{rcll}
   \Psi_{\frac{17+5}{2}}=\Psi_{11} & \equiv &
   -6z_4(z_4^6-7 z_4^3z_6^2+6z_6^4)& \text{mod}\ \  17\\
    \ & \ & \ &\ \\
   \  & \equiv
   &  -6z_4(z_4^3-z_6^2)(z_4^3-6z_6^2)&  \text{mod}\ \ 17,\\
   \ & \ & \ &\ \\
   \Delta & \equiv  & 4(z_4^3-6z_6^2)& \text{mod}\ \ 17,\\
    \ & \ & \ &\ \\
   H=H_{17} & \equiv  & 2z_4(z_4^3-z_6^2)&  \text{mod}\ \ 17, \end{array}$$
   hence
    \begin{equation}
    \label{bb17}
\Psi_{11} \equiv  10\Delta H_{17}\ \ \   \text{mod}\ \ 17.
   \end{equation}
    \medskip

  \noindent Congruences \ref{bb11}, \ref{bb13}, \ref{bb17} plus our Theorem \ref{godem1}, imply our Corollary \ref{godem4}. 
  
   It would be interesting to understand if  congruences of the type 
    \ref{bb11}, \ref{bb13}, \ref{bb17}
    hold for (almost) all primes.
Note that these congruences are not simply  consequences of equalities in ${\mathbb Q}[z_4,z_6]$: indeed 
   $\Psi_9$ and $\Delta H_{13}$ are ${\mathbb Q}$-linearly independent,
   for there are sign changes in the coefficients of $\Psi_9$ but no sign changes
   in the coefficients of $\Delta H_{13}$. Similarly a direct calculation shows that   $\Psi_8$ and $\Delta H_{11}$ are    ${\mathbb Q}$-linearly independent.  So whatever lies behind congruences \ref{bb11}, \ref{bb13}, \ref{bb17} seems to be of a more mysterious nature.  \end{remark}
  
  The next Proposition shows that $\Psi_{\frac{p+5}{2}}\not \equiv 0$ mod $p$ for 
  $p\not\equiv 11$ mod $12$; the case $p\equiv 11$ mod $12$ seems to be significantly trickier  to handle.

   \begin{proposition}\label{golem}
   Let $p\geq 13$.
   
   1) If $p\equiv 1$ mod $3$ the following holds:
   $$\Psi_{\frac{p+5}{2}}(0,1)\not\equiv 0\ \ \ \text{mod}\ \ p.$$
   
   2) If $p\equiv 1$ mod $4$ the following holds:
   $$\Psi_{\frac{p+5}{2}}(1,0)\not\equiv 0\ \ \ \text{mod}\ \ p.$$
   
   \end{proposition}
   
   {\it Proof}. Assume  $p$ is arbitrary.
   By our recurrence relation \ref{tututu} for the $\psi_n$'s we have
   $$\Psi_n(0,1)\equiv  \frac{2n-7}{2n+2}\cdot \frac{2n-9}{2n}\cdot \Psi_{n-3}(0,1) \ \ \ \text{mod}\ \ p$$
   for $5\leq n\leq p-2$. Note also that for $n\leq \frac{p+5}{2}$ we have 
   $$2n-9\leq 2n-7\leq p-2$$
   hence, for $n$ in this range we have that $\Psi_n(0,1)\equiv 0$ mod $p$ if and only if $\Psi_{n-3}(0,1)\equiv 0$ mod $p$. 
 Assume  $p\equiv 1$ mod $3$ and  set $p=6N+1$.
   Since
   $$\Psi_6(0,1)\equiv \frac{15}{2^7\times 7}\not\equiv 0\ \ \text{mod}\ \ p,$$
we get
      $$\Psi_{\frac{p+5}{2}}(0,1)=\Psi_{3N+3}(0,1)\not\equiv 0\ \ \text{mod}\ \ p,$$
      which ends our proof in case 1).
      
      Assume, again, that $p$ is arbitrary.
       For $5\leq 2n-1\leq p-2$ our recurrence relation
      \ref{tututu} for the $\psi_n$'s gives
      $$\Psi_{2n-1}(1,0)\equiv -\frac{4n-9}{4n}\cdot \frac{4n-5}{4n-2}\cdot  \Psi_{2n-3}(1,0)\ \ \ \text{mod}\ \ p.$$
      For $2n-1\leq \frac{p+5}{2}$ we have 
      $$4n-9\leq 4n-5\leq p+2.$$
      Now the only positive integer $\leq p+2$ divisible by $p$ is $p$ itself and none of $4n-9$ or $4n-5$ can be equal to $p$ because $p\equiv 1$ mod $4$. We get that $\Psi_{2n-1}(1,0)\equiv 0$ mod $p$ if and only if $\Psi_{2n-3}(1,0)\not\equiv 0$ mod $p$. 
      Assume $p\equiv 1$ mod $4$ and set $p=4N+1$. Since
      $$\Psi_5(1,0)\equiv -\frac{7}{2^8\times 5}\not\equiv 0\ \ \ \text{mod}\ \ p,$$
      we get
      $$\Psi_{\frac{p+5}{2}}(1,0)=\Psi_{2N+3}(1,0)\not\equiv 0\ \ \ \text{mod}\ \ p,$$
      which ends our proof in this case as well.
   \qed

   \medskip
      
   {\it Proof of Theorem \ref{godem5} and  Proposition \ref{supplement}}.
   These  follow directly from
  Proposition \ref{golem}, 
   Theorem \ref{godem1}, and Proposition \ref{cal}. \qed
   
   \medskip
   
   {\it Proof of Proposition \ref{supplement2}.}
   Assume in what follows that $p, a, b$ are arbitrary subject to the conditions:
 $$p\equiv 1\ \ \  \text{mod}\ \  3, \ \ \ a\equiv 0\ \ \text{mod}\ \ p,\ \ b\not\equiv 0\ \ \text{mod}\ \ p.$$
 Set $p=3N+1$; then $N$ is, of course, even.
  Note that 
 $$\Delta(a,b)\equiv 27b^2 \not\equiv 0\ \ \ \text{mod}\ \ \ p,$$
and 
$$H(a,b)\equiv \left(\begin{array}{c}
\frac{p-1}{2}\\ 
\ \\
\frac{p-1}{3}\end{array}\right)b^{\frac{p-1}{6}}\not\equiv 0\ \ \ \text{mod}\ \ \ p.$$
Since $\Sigma=z_6\Psi H$ with $\Psi=\Psi_{\frac{p+5}{2}}$ and $\Psi(a,b)\not\equiv 0$ mod $p$ (cf. Proposition \ref{golem}) it follows that $a,b$ is $\Sigma$-non-singular. Consider again the congruences \ref{planeplan} as a linear system ${\mathfrak R}$ with ``unknowns" $v_0,\theta$.
By Proposition \ref{cal} this system is uniquely solvable mod $p$. We will prove in what follows the following:

\medskip

{\it Claim. 
There is a universal $\beta\in \bZ_{(p)}$ such that 
$v_0=0,\ \theta=-\frac{\d b}{6b^p}-\beta$
is a solution to the system ${\mathfrak R}$ mod $p$.}

\medskip

 Assuming this is the case it follows that, with $\Theta$ as in the proof of Theorem \ref{godem1}, we have
$$\Theta(0,z_6,z'_4,z'_6)\equiv -\frac{z'_6}{6z_6^p}-\beta\ \ \ \text{mod}\ \ p,$$
and  Proposition \ref{supplement2} follows.

We now prove our Claim.

Under our hypothesis $a\equiv 0$ mod $p$ the system ${\mathfrak S}_{\frac{p+7}{2}}$, i.e. the one defined by \ref{louvre} for $s\leq \frac{p+7}{2}$, becomes:
\begin{equation}
 \label{louvre2}
 sb^pv_s\equiv 
 \left(\frac{9}{2}-s\right)v_{s-3}+c_s+d_s+e_s+f_s\ \ \ \text{mod}\ \ \ p,\ \ \ 1\leq s\leq \frac{p+7}{2},
 \end{equation}
 where 
 $$c_4\equiv \frac{3\theta}{2},\ \ e_1\equiv \frac{\d b}{2},\ \ f_2\equiv \frac{\d a}{2}\ \ \text{mod}\ \ p,$$
 all the other $c_s, e_s, f_s\equiv 0$ mod $p$, and
 $d_s=D_s(0,b)$ with $D_s$ universal functions  in $M_H$ of weight $(8-2s)p$
 that vanish  for $s\geq 5$. So we have
 $$d_1\equiv \beta_1b^p,\ \ d_2\equiv 0,\ \ d_3\equiv 0,\ \ \ d_4\equiv \beta_4\ \ \text{mod}\ \ p$$
 where $\beta_1, \beta_4\in \bZ_{(p)}$ are some universal constants.
 Set
 \begin{equation}\label{thecond}
 \beta:=\frac{1}{3}\beta_1+\frac{2}{3}\beta_4,\ \ \ \ 
 \theta\equiv - \frac{\d b}{6b^p}-\beta\ \ \text{mod}\ \ p.\end{equation}
  The system ${\mathfrak S}_{\frac{p+7}{2}}$ is then equivalent to:
 
 \medskip

 \begin{equation}
 \label{solsol}
 \begin{array}{lllllll}
 b^pv_1 & \equiv & \frac{\d b}{2}+\beta_1b^p, & \ & 2 b^p v_2 & \equiv &  \frac{\d a}{2},\\
    \ & \ & \ & \ & \ & \  & \ \\
3b^p v_3 & \equiv & \frac{3}{2}v_0, & \ & 4b^p v_4 & \equiv & \frac{\d b}{2^2b^p}+\frac{3\theta}{2}+\frac{1}{2}\beta_1+\beta_4,\\
 \ & \ & \ & \ & \ & \  & \ \\
  5b^p v_5 & \equiv & -\frac{1}{2}v_2, &\ & 6b^p v_6 & \equiv & -\frac{3}{2} v_3,\\
   \ & \ & \ & \ & \ & \  & \ \\
7b^p v_7 & \equiv & -\frac{5}{2} v_4, &\ & 8b^p v_8 & \equiv & - \frac{7}{2} v_5,\\
  \ & \ & \ & \ & \ & \  & \ \\
  9b^p v_9 & \equiv & -\frac{9}{2}v_6, &\ & \text{etc}. & \  &
  \end{array}
 \end{equation} 
This system has a unique (mod $p$) solution
 $v_0,v_1,v_2,...,v_{\frac{p+7}{2}}$ 
 with $v_0=0$. This solution has the property that 
  $v_{3k}=0$ for  $0\leq 3k\leq \frac{p+7}{2}$ and $v_{3k+1}=0$ for  $1\leq 3k+1\leq \frac{p+7}{2}$. In particular
 $$v_{\frac{p+5}{2}}=v_{3(\frac{N}{2}+1)}\equiv 0\ \ \text{mod}\ \ p,$$
 $$v_{\frac{p+7}{2}}=v_{3(\frac{N}{2}+1)+1}\equiv 0\ \ \ \text{mod}\ \ p.$$
 So the system ${\mathfrak R}$ is satisfied by $v_0=0$ and $\theta=-\frac{\d b}{6b^p}-\beta$;
 this proves the Claim above and we are done.\qed

 \medskip
 
 {\it Proof of Proposition \ref{5trays}}.
 We assume in what follows that
 $$p\equiv 1\ \ \  \text{mod}\ \  4, \ \ \ a\not\equiv 0\ \ \text{mod}\ \ p,\ \ b\equiv 0\ \ \text{mod}\ \ p.$$
Set $p=4N+1$.
 Note that $a,b$ is then an ordinary pair because
 $$\Delta(a,b)\equiv 4a^3 \not\equiv 0\ \ \ \text{mod}\ \ \ p,$$
and 
$$H(a,b)\equiv \left(\begin{array}{c}
\frac{p-1}{2}\\ 
\ \\
\frac{p-1}{4}\end{array}\right)a^{\frac{p-1}{4}}\not\equiv 0\ \ \ \text{mod}\ \ \ p.$$
Assume from now on $\Sigma=H$.
The system ${\mathfrak S}$ deduced from \ref{louvre} becomes:
\begin{equation}
 \label{louvre3}
 \left(s-\frac{3}{2}\right) a^pv_{s-1} \equiv 
 \left(\frac{9}{2}-s\right)v_{s-3}+c_s+d_s+e_s+f_s\ \ \ \text{mod}\ \ \ p,\ \ \ s\geq 1,
 \end{equation}
 where 
 $$c_2\equiv \frac{a^p\theta}{2},\ \ c_4\equiv \frac{3\theta}{2},\ \ 
 e_1\equiv \frac{\d b}{2},\ \ f_2\equiv \frac{\d a}{2},\ \ 
 d_2\equiv \alpha_2a^p,\ \ d_4\equiv \alpha_4\ \ \text{mod}\ \ p,$$
  where $\alpha_2, \alpha_4\in \bZ_{(p)}$ are some universal constants and all the rest of $c_s,d_s,e_s,f_s$ are $\equiv 0$ mod $p$. Explicitly, the truncation ${\mathfrak S}_5$ of ${\mathfrak S}$ reads:
  
  \medskip

   \begin{equation}
 \label{solsolsol}
 \begin{array}{rllllll}
 -\frac{1}{2}a^pv_0 & \equiv & \frac{\d b}{2}, & \ &
   \frac{1}{2} a^p v_1 & \equiv &  \frac{a^p\theta}{2}+\frac{\d a}{2}+\alpha_2 a^p,\\
    \ & \ & \ & \ & \ & \ & \ \\
\frac{3}{2}a^p v_2 & \equiv & \frac{3}{2}v_0, & \ & \frac{5}{2}a^p v_3 & \equiv & \frac{\d a}{2a^p}+2\theta+\alpha_2+\alpha_4,\\
  \ & \ & \ & \ & \ & \ & \ \\
  \frac{7}{2} a^p v_4 & \equiv & -\frac{1}{2}v_2, & \ & \frac{9}{2}a^p v_5 & \equiv & -\frac{3}{2} v_3.
  \end{array}
 \end{equation} 
Set
 \begin{equation}
 \label{thecond4}
 \alpha:=\frac{1}{2}\alpha_2+\frac{1}{2}\alpha_4,\ \ \ 
 \theta\equiv - \frac{\d a}{4a^p}-\alpha\ \ \text{mod}\ \ p.\end{equation}
 By Proposition \ref{operastabbing} we will be done if we 
 can construct a sequence of elements $v_0,v_1,v_2,...\in R$, almost all of which are $0$, satisfying the congruences \ref{louvre3}.
 Define $v_0,v_1,v_2,v_3,v_4,v_5$ via the conditions \ref{solsolsol}. 
 In particular, by \ref{thecond4},  $v_3\equiv 0$ mod $p$. Then take $v_{2k+1}\equiv 0$ for all 
 $k\geq 2$; so \ref{louvre3} holds for all odd $s$. Now define
 $v_{2k}$ via 
 \begin{equation}
 \label{pattern2k}
  v_{2k}\equiv (-1)^{k+1} \frac{5\times 9 \times 13\times  ...
 \times (4k-7)}{7\times 11\times 15 \times ... \times (4k-1)}\cdot \frac{v_2}{a^{(k-3)p}},
 \end{equation}
  for $3\leq k\leq k_1:=N+2$. This is allowable because of the following:
  
  \medskip
  
  {\it Claim}. For $k\leq k_1$ we have    $4k-1\not\equiv 0$ mod $p$.
  
  \medskip
  
   To check the Claim note that
 for $k\leq k_1$, we have
  $4k-1\leq p+6$ so the only value of $4k-1$ in this range divisible by $p$ is $p$ itself; but the equality $4k-1=p$ is impossible, which ends the proof of the Claim.
  
   Now note that $v_{2k_1}\equiv 0$ mod $p$ because $4k_1-7=p$. So we may take $v_{2k}=0$ for all $k\geq k_1$ and the congruences \ref{louvre3} will be satisfied for all even $s$, which concludes our proof.\qed
   
   \medskip   
      
{\it Proof of Remark \ref{cryterion}}.
We will check assertion 2; assertion 1 can be proved similarly.
 For   any $\Sigma$-non-singular pair $a,b\in R$, and  any $c\in R^{\times}$, we have
the following congruences mod $p^2$:
$$\begin{array}{rcl}
F(c^4a, c^6b) & \equiv &
c^k\Gamma_k(a,b)\\
\ & \ & \ \\
\ & \ & +
pc^{k-4p}\Gamma^*_{k-4p}(a,b)\d(c^4a)\\
\ & \ & \\
\ & \ & +pc^{k-6p}\Gamma^*_{k-6p}(a,b)\d(c^6b)\\
\ & \ & \ \\
\ & \equiv & c^k\Gamma_k(a,b)\\
\ & \ & \ \\
\ & \ & +
c^{k-4p}\Gamma^*_{k-4p}(a,b)((c^p+p\d c)^4(a^p+p\d a)-c^{4p}a^p)\\
\ & \ & \ \\
\ & \ & +c^{k-6p}\Gamma^*_{k-6p}(a,b)(
(c^p+p\d c)^6(b^p+p\d b)-c^{6p}b^p
)\\
\ & \ & \ \\
\ & \equiv & c^k\Gamma_k(a,b)\\
\ & \ & \ \\
\ & \ & +
pc^{k-4p}\Gamma^*_{k-4p}(a,b)(4c^{3p}a^p\d c+c^{4p}\d a)\\
\ & \ & \ \\
\ & \ & +pc^{k-6p}\Gamma^*_{k-6p}(a,b)(6c^{5p}b^p\d c+c^{6p}\d b)\\
\ & \ & \ \\
\ & \equiv & c^k\Gamma_k(a,b)\\
\ & \ & \ \\
\ & \ & +p(4a^p\Gamma^*_{k-4p}(a,b)+6b^p\Gamma^*_{k-6p}(a,b))
c^{k-p}\d c\\
\  & \ & \\
\ & \ & + pc^k\Gamma^*_{k-4p}(a,b)\d a+
 pc^k\Gamma^*_{k-6p}(a,b)\d b \ \ \ \text{mod}\ \ \ p^2,\end{array}$$
 On the other hand we have:
$$
\begin{array}{rcl}
c^{k+p-\phi}F(a,b) & \equiv &
 c^k(1-p\frac{\d c}{c^p})(\Gamma_k(a,b)+
 p\Gamma^*_{k-4p}(a,b)\d a+p\Gamma^*_{k-6p}(a,b)\d b)\\
\ & \ & \ \\
\ & \equiv & c^k\Gamma_k(a,b)-p\Gamma_k(a,b)c^{k-p}\d c\\
\ & \ & \ \\
\ & \ & +
pc^k\Gamma^*_{k-4p}(a,b)\d a+pc^k\Gamma^*_{k-6p}(a,b)\d b.\ \ \ \text{mod}\ \ \ p^2.
\end{array}
$$

\medskip

\noindent Our Remark follows. \qed

\end{document}